\pdfoutput=1
\documentclass[11pt]{amsart}
\usepackage{graphicx}
\usepackage{amsmath,amssymb,amsfonts,amscd,latexsym}

\newcommand{\labell}[1] {\label{#1}}

\newcommand{\1}{{{\mathchoice {\rm 1\mskip-4mu l} {\rm 1\mskip-4mu l}
{\rm 1\mskip-4.5mu l} {\rm 1\mskip-5mu l}}}}

\newcommand{\rank}{{\rm rank}}
\newcommand{\Tg}{{\Tilde g}}
\newcommand{\TX}{{\Tilde X}}
\newcommand{\Qq}{{\mathcal Q}}
\newcommand{\PD}{{\rm PD}}

\newcommand{\ovv}{{\overline v}}
\newcommand{\ovV}{{\overline V}}

\newcommand{\Trho}{{\Tilde \rho}}

\newcommand{\Dd}{{\mathcal D}}

\newcommand{\Ver}{{\rm Vert}}

\newcommand{\less}{{\smallsetminus}}

\newcommand{\bla}{{\bigl\langle}}
\newcommand{\bra}{{\bigl\rangle}}

\newcommand{\Tga}{{\Tilde\ga}}
\newcommand{\THam}{{\Tilde\Ham}}

\newcommand{\al}{{\alpha}}

\newcommand{\be}{{\beta}}

\newcommand{\Om}{{\Omega}}
\newcommand{\om}{{\omega}}
\newcommand{\eps}{{\varepsilon}}
\newcommand{\de}{{\delta}}

\newcommand{\ga}{{\gamma}}
\newcommand{\Ga}{{\Gamma}}
\newcommand{\io}{{\iota}}
\newcommand{\ka}{{\kappa}}
\newcommand{\la}{{\lambda}}
\newcommand{\La}{{\Lambda}}
\newcommand{\si}{{\sigma}}

\newcommand{\Mm}{{\mathcal M}}
\newcommand{\Ss}{{\mathcal S}}
\newcommand{\Ll}{{\mathcal L}}
\newcommand{\oMm}{{\overline {\Mm}}}
\newcommand{\ov}{\overline}
\newcommand{\oc}{\ov {c}}

\newcommand{\id}{{\rm id}}

\renewcommand{\Tilde}{\widetilde}

\newcommand{\Tphi}{{\Tilde \phi}}
\newcommand{\Tpsi}{{\Tilde \psi}}

\newcommand{\TP}{{\Tilde P}}

\newcommand{\Tom}{{\Tilde \om}}
\newcommand{\TOm}{{\Tilde \Om}}
\newcommand{\Tsi}{{\Tilde \si}}

\newcommand{\Q}{{\mathbb Q}}
\newcommand{\R}{{\mathbb R}}
\newcommand{\C}{{\mathbb C}}
\newcommand{\Crit}{{\rm Crit}}

\newcommand{\Z}{{\mathbb Z}}

\newcommand{\Ham}{{\rm Ham}}

\newcommand{\Vol}{{\rm Vol}}

\newcommand{\Aa}{{\mathcal A}}

\newcommand{\SSS}{{\smallskip}}
\newcommand{\QED}{{\hfill $\Box$\MS}}

\newtheorem{theorem}{Theorem}[section]
\newtheorem{thm}[theorem]{Theorem}

\newtheorem{cor}[theorem]{Corollary}
\newtheorem{lemma}[theorem]{Lemma}

\newtheorem{prop}[theorem]{Proposition}

\newtheorem{defn}[theorem]{Definition}

\newtheorem{rmk}[theorem]{Remark}

\numberwithin{figure}{section}
\numberwithin{equation}{section}
\numberwithin{table}{section}

\newcommand{\MS}{{\medskip}}

\newcommand{\NI}{{\noindent}}

\begin{document}

 \title{Monodromy in Hamiltonian Floer theory}
 \author{Dusa McDuff}\thanks{partially supported by NSF grant DMS 0604769.}
\address{Department of Mathematics,
Barnard College, Columbia University, New York, NY 10027-6598, USA.}
\email{dmcduff@barnard.edu}
\date{26 December 2007, last revised  February 15, 2010}

\begin{abstract}  Schwarz showed  that when
a closed symplectic manifold $(M,\om)$ is symplectically aspherical (i.e. the symplectic form and the first Chern class vanish  on $\pi_2(M)$) then
the spectral invariants, which are initially defined on the universal cover
of the Hamiltonian group, descend  to the Hamiltonian group 
Ham$(M,\om)$.
 In this note we describe less stringent conditions
  on the Chern class and quantum  homology of $M$ under which the
  (asymptotic) spectral invariants descend to Ham$(M,\om)$. For example, they  
  descend if the quantum multiplication of
   $M$ is undeformed and $H_2(M)$ has rank $>1$, or if the minimal Chern number is at least $n+1$ (where $\dim M=2n$) and the even cohomology
  of $M$ is generated by divisors. 
  The proofs are based on certain calculations of genus zero
   Gromov--Witten invariants.
   As an application, we show that the Hamiltonian group 
   of the one point blow up of $T^4$ 
   admits a Calabi quasimorphism.  Moreover, whenever the 
   (asymptotic) spectral invariants descend it is easy to see that 
   Ham$(M,\om)$ has infinite diameter in the Hofer norm. Hence our results 
   establish the infinite diameter of Ham in many new cases. We also 
   show that the area pseudonorm --- a geometric version  of the Hofer norm --- is nontrivial on the
   (compactly supported)  Hamiltonian group  for all noncompact manifolds
   as well as for a large class of closed manifolds.
    \end{abstract}
 
\maketitle
\tableofcontents

\section{Introduction}\labell{s:intro}

Let $(M,\om)$ be a closed symplectic manifold.  
Denote by $\Ham: = \Ham(M,\om)$ its group of Hamiltonian symplectomorphisms and by  
$\THam$ the universal cover of $\Ham$.  Each path $\{\phi_t\}_{0\le t\le 1}$ in $\Ham$ is the flow of some time dependent Hamiltonian $H_t$ and, following Hofer \cite{HOF1}, we define its  length $\Ll(\{\phi_t\})$ to be:
$$
\Ll(\{\phi_t\}) = \int_0^1 \Bigl(
\max_{x\in M}H_t(x) - \min_{x\in M}H_t(x)\Bigr)dt.
$$
The Hofer (pseudo)norm $\|\Tphi\|$  of an element $\Tphi = (\phi, \{\phi_t\})$ in 
the universal cover $\THam$
of $\Ham$ is then defined to be the infimum of the lengths of the paths 
from the identity element $id$
to $\phi$ that are homotopic to $\{\phi_t\}$.   Similarly  we define the norm $\|\phi\|$ of an element in $\Ham$ to be the infimum of the lengths of {\it all} paths from $id$ to $\phi$.  It is easy to see that
$\|\phi\|$ is conjugation invariant and satisfies $\|\phi\psi\|\le \|\phi\|+\|\psi\|$, but harder to see that it is nondegenerate, i.e. $\|\phi\|=0$ iff $\phi=id$.  (This was proved for compactly supported symplectomorphisms of $\R^{2n}$ by Hofer \cite{HOF1} and for general
 $M$ by Lalonde--McDuff \cite{LMe}.)
It is unknown whether $\|\cdot\|$ is always nondegenerate (and hence a norm) on $\THam$ since it may vanish on some elements of the subgroup $\pi_1(\Ham)$.
On the other hand, there is no known counterexample; for some results in the positive direction see Remark~\ref{rmk:cN} below.

The question of whether $\|\cdot\|$ is uniformly bounded makes sense 
even if $\|\cdot\|$ is just a pseudonorm.   If is it unbounded 
on $\THam$ or $\Ham$   we shall say that this group
 has {\it infinite (Hofer) diameter.}
 Ostrover~\cite{Os} showed
 that $\THam$ always  has infinite Hofer diameter
  (we sketch the proof below), while the corresponding 
  result is unknown for $\Ham$
in many cases.  For example, it is shown in \cite{Mcsn} that $\Ham$ has infinite diameter when $M$ is a \lq\lq small" 
blow up of $\C P^2$ but 
it is unknown whether this remains the case when $M$ is monotone
(i.e. the exceptional divisor is precisely one third the size of the line) or is a still
bigger blow up.  

However, if both $[\om]$ and 
$c_1(M)$ vanish on $\pi_2(M)$
then Schwarz \cite{Sch} showed that $\Ham$ does have infinite diameter.  To prove this he established that for such $M$ each element
$\Tphi\in \THam$ has  a set of so-called
{\it spectral invariants} 
$$
\{c(a,\Tphi)\,|\,a\in QH_*(M), \;a\ne 0\}\subset \R.
$$
Later, Oh~\cite{Oh,Oh1} and Usher~\cite{U} 
showed that the numbers $c(a,\Tphi)$ are  
well defined on $\THam$  for all symplectic manifolds.
It follows easily from their properties 
(explained in \S\ref{s:si} below) that $\THam$ has infinite Hofer diameter.
However, 
 they do not in general descend to well defined functions on $\Ham$; in other words it may not be true that   $c(a,\Tphi)=c(a,\Tpsi)$ whenever $\Tphi,\Tpsi$ project to the same element of $\Ham$.
Schwarz showed  that when both $[\om]$ and 
$c_1(M)$ vanish on $\pi_2(M)$ the invariants do descend to $\Ham$. As we explain below, it is then an easy consequence of Ostrover's construction
 that $\Ham$ has infinite diameter.

Another case in which $\Ham$ was known to have infinite diameter is
 $S^2$.  The original proof in Polterovich \cite{P} initially appears somewhat different in spirit  from the approach 
  presented here, but the  arguments  in Entov--Polterovich~\cite{EP1}
  using quasimorphisms
bring  this result also within the current framework; 
cf. case (ii) of Theorem \ref{thm:main2}.

In fact there are two questions one can ask here.  Do the spectral invariants themselves descend, or is it only their asymptotic versions 
that descend?  Here, following Entov--Polterovich~\cite{EP1}, 
we define the {\it  asymptotic spectral invariants} $\oc(a,\Tphi)$ for nonzero $a\in QH_*(M)$ by setting
$$
\oc(a,\Tphi): = \liminf_{k\to \infty} \frac{c(a,\Tphi\,^k)}{k},\quad\mbox{ for all }\;\;\Tphi\in \THam.
$$
Ostrover's construction again implies that  $\Ham$ has
 infinite diameter whenever an asymptotic spectral number
descends to $\Ham$; see Lemma~\ref{le:inf}.  

We shall extend Schwarz's result in two directions, imposing conditions
either on $\om$ via the Gromov--Witten invariants or on $c_1$. 
Recall that  the quantum product on $H_*(M)$ is defined by  
the $3$-point
genus zero Gromov--Witten invariants 
$$
\bla a_1,a_2,a_3\bra^M_\be, \qquad \be\in H_2(M), a_i\in H_*(M),
$$ 
and reduces to the usual intersection product if these invariants vanish
whenever $\be\ne 0$. In the latter case we shall say that
the  quantum product (or simply $QH_*(M)$) is {\it
undeformed}.
The condition $[\om]|_{\pi_2(M)}=0$  is much stronger; in this case  there are no $J$-holomorphic spheres at all, and so the quantum product
on $M$ is of necessity undeformed.  If $[\om]|_{\pi_2(M)}=0$,
Schwarz's argument easily extends to show that the spectral invariants descend; see Proposition~\ref{prop:sch}(i). However,  its generalization 
in  Theorem~\ref{thm:main1} below concerns the asymptotic invariants.

We shall denote by $N$
 the {\it  minimal Chern number} of $(M,\om)$, i.e. 
the smallest positive value of  $c_1(M)$ on $\pi_2(M)$.
If $c_1|_{\pi_2(M)}=0$  then we set $N: = \infty$. 
We always denote $\dim\, M = 2n$.
We shall say that $(M,\om)$ is {\it spherically monotone} if 
there is  $\ka>0$ such that
$c_1|_{\pi_2(M)}=\ka\, \om|_{\pi_2(M)}$  and is {\it negatively 
monotone} if 
$c_1=\ka\, [\om]$  on $\pi_2(M)$  for some $\ka< 0$.
Recall also that
$(M,\om)$  is said to be {\it (symplecically) uniruled} if some 
genus zero  Gromov--Witten invariant
of the form $\bla pt,a_2,\dots,a_m\bra_\be$, $\be\ne 0$, does not vanish.  It is called {\it strongly uniruled} if this happens for $m=3$.

\begin{thm}\labell{thm:main1}  Suppose that  $(M,\om)$ is a closed symplectic manifold such that
$QH_*(M)$ is undeformed.  Then the asymptotic spectral invariants
descend to $\Ham$ except possibly if 
the following three additional conditions all hold: 
$\rank\, H_2(M) = 1,  N\le n,$
and $(M,\om)$ is spherically monotone.
\end{thm}

\begin{rmk}\rm (i)
The exceptional case does not occur if $M$ has dimension $4$.
For if $QH_*(M)$ is undeformed then $(M,\om)$ is minimal. (The class of an exceptional sphere always has nontrivial Gromov--Witten invariant; see \cite{Mcu}.) Moreover it follows from the results of Taubes--Li--Liu that it cannot be  spherically monotone,
for if it were it would be uniruled, in particular the quantum product would be deformed. 
Similarly,  the exceptional case does not occur for 
smooth projective varieties of any dimension.  For $N\ne \infty$ implies that
at least some nonzero elements in $H_2(M)$ are represented by spheres.
Hence they must all be (since $\rank\, H_2(M) = 1$).  Therefore the conditions imply that $M$ is Fano, and hence, by an argument of
Kollar--Ruan that is explained in \cite{HLR}, also symplectically uniruled.
   It is unknown whether
 every spherically monotone symplectic manifold is uniruled. 
  If this were true then there would be no exceptions at all.\MS
  
  \NI (ii)  Every minimal $4$-manifold that is not rational or ruled (such as a $K3$ surface or a simply connected surface of general type) has vanishing genus zero Gromov--Witten invariants and so is covered by this theorem.
\end{rmk}

The next result explains what can be proved under various conditions on 
the minimal Chern number $N$. 
We denote the  even degree homology of $M$ by  $H_{ev}(M)$.

\begin{thm}\label{thm:main2} Let $(M,\om)$ be a closed symplectic $2n$-dimensional manifold with minimal Chern number $N$.  Suppose further that either $H_{ev}(M)$ is generated as a ring by 
the divisors $H_{2n-2}(M)$ or that $QH_*(M)$ is undeformed.  Then: \SSS

\NI
{\rm (i)}  If $N\ge n+1$, the spectral invariants are well 
defined on $\Ham$ except possibly if $N\le 2n$ and
$(M,\om)$ is strongly uniruled.
\SSS

\NI {\rm (ii)}  If $n+1\le N \le 2n$ the asymptotic 
spectral invariants are well defined on $\Ham$.\SSS

\NI {\rm (iii)}  The conclusion in {\rm (ii)} still holds if $N=n$ except possibly if
$(M,\om)$ is strongly uniruled  or if $\rank\, H_2(M)=1$.\SSS

\NI {\rm (iv)}  The conclusion in {\rm (ii)} also holds 
when $(M,\om)$ is negatively monotone, independently of the values of $n,N$.
\end{thm}

For example, if $M$ is a $6$-dimensional K\" ahler manifold then
$H_{ev}(M)$ is generated as a ring by $H_{4}(M)$ because 
$\wedge [\om]: H^2(M)\to H^4(M)$ is an isomorphism.  Hence all
Calabi--Yau $3$-folds satisfy the conditions of this theorem.

\begin{rmk}\rm \labell{rmk:main} (i) In Theorem~\ref{thm:main2} the conditions in the second sentence
 may be replaced by the weaker but somewhat technical
  condition (D); cf. Definition~\ref{def:D}. \SSS

\NI (ii) Theorem~\ref{thm:main2} is
sharp.  To see that the non-uniruled
 hypothesis in (i)  is necessary,  observe that by 
 Entov--Polterovich~\cite{EP1} the spectral invariants do not descend
 for $M=\C P^n$ although the asymptotic ones do.  
 This condition  is also needed in (iii).  For example,
consider  $M=S^2\times S^2$ which has $N=n=2$
   and is strongly uniruled.
  Ostrover~\cite{Os2} showed that 
 the asymptotic spectral invariants descend   if and only if
$(M,\om)$ is  monotone, i.e. the two $2$-spheres have equal area. 
Further, the results do not extend to smaller $N$.
 Proposition~\ref{prop:bl} below gives many examples 
of manifolds with $N=n-1$ that satisfy the other cohomological conditions but
are such that the asymptotic invariants do not descend.
\SSS

\NI (iii)   It is not clear what happens when
$N=n$ but 
$\rank\,H_2(M) =1$.  If $(M,\om)$ is negatively monotone, then 
by Proposition~\ref{prop:sch} the asymptotic spectral invariants always descend (without any condition on $QH_*(M)$), but it is not clear what happens 
 in the positive case. (Cf. the similar missing 
 case in Theorem~\ref{thm:main1}.)
If $(M,\om)$  were also projective then
$(M,\om)$ would be uniruled and one would not 
expect the invariants to descend but in the general case considered here all we can say is that our methods fail.
The relevant part of the proof 
of Proposition~\ref{prop:sch} fails  for $N=n$ and $\ka>0$, while the argument in 
Lemmas~\ref{le:hard} and \ref{le:hard1} definitely needs  $\rank\,H_2(M) >1$. \SSS
 
 \NI (iv) The example of $S^2\times S^2$ in Remark~\ref{rmk:main}(ii) 
above suggests that perhaps
 the  asymptotic spectral invariants descend for all 
  monotone manifolds.
 But this is not true.  Consider, for example the monotone one point blow up of $\C P^2$ with its obvious $T^2$ action. It is easy to see 
 that there are  circles in $T^2$ that represent elements 
 $\ga\in \pi_1(\Ham)$ for 
 which $\oc(\1,\ga)\ne 0$;
 see  \cite{Mcsn,Os2}.\SSS
 \end{rmk}

\begin{cor} \labell{cor:main}  If $(M,\om)$ satisfies any of the conditions in
Theorems~\ref{thm:main1} and~\ref{thm:main2}, then $\Ham$ has infinite Hofer diameter. 
\end{cor}
\begin{proof}  This holds by Lemma~\ref{le:inf}.
\end{proof}

Of course, one expects $\Ham$ always to have infinite  Hofer diameter, but this question seems out of reach with current techniques.   However there are other ways to tackle this question. For example, in \cite{Mcsn}  we show that a small blow up of $\C P^2$ has infinite Hofer diameter even though the spectral invariants do not descend by using an argument based on the asymmetry of the spectral invariants, i.e. the fact that the function $\ovV$ of Remark~\ref{rmk:spectv} does  not vanish on $\pi_1(\Ham)$. Also if 
$\pi_1(M)$ is infinite, one can sometimes use the energy--capacity inequality  as in Lalonde--McDuff \cite{LM}.

Another related problem is the question of when the 
{\it area pseudonorm}
$\rho^+ + \rho^-:\Ham\to \R$ defined in \cite{Mcv}
 is nonzero. 
  Here 
$$
\rho^+(\phi): = \inf \;\int_0^1\max_{x\in M}H_t \,dt,
$$
where the infimum is taken over all mean normalized\footnote{
i.e $\int_M H_t\om^n=0$ for all $t$.  Also, this discussion of one sided norms is the one place in this paper
where the choice of signs is crucial. In order to be consistent with \cite{Mcv}  we shall define the flow of $H_t$ to be generated by the vector $X_{H_t}$ satisfying  $\om(X_H,\cdot) = -dH_t$.}
Hamiltonians with time $1$ map $\phi$.  Further
$\rho^-(\phi): = \rho^+(\phi^{-1})$. It is easy to see that
$\rho^+ + \rho^-$ is a conjugation invariant pseudonorm
on $\Ham$.  Therefore, because $\Ham$ is simple, $\rho^+ + \rho^-$
 is either
 identically zero or is
nondegenerate and hence a norm. 
  The difficulty in dealing with it is that one may need to use 
different lifts $\Tphi$ of $\phi$ to $\THam$ to calculate
 $\rho^+$ and $\rho^-$.  However it has a very natural geometric interpretation. For by \cite[Prop.~1.12]{Mcv}  
 $$
 \rho^+(\phi) + \rho^-(\phi)=\inf\,\bigl(\Vol(P,\Om)/\Vol(M,\om)\bigr)
 $$
 where the infimum is taken over all  Hamiltonian fibrations
 $(M,\om)\to (P,\Om)\to S^2$  with 
 monodromy $\phi$ around some embedded loop in the base. This ratio 
is called the {\it area} of the fibration $P\to S^2$ 
since for product fibrations it would be the area of the base.

 Since the area pseudonorm
  $ \rho^+ + \rho^-$ is never larger than the Hofer norm
 the next result also implies Corollary~\ref{cor:main}.
 
 \begin{cor}\labell{cor:rho+}   
 If $(M,\om)$ satisfies 
 any of the conditions in
Theorems~\ref{thm:main1} and~\ref{thm:main2}, then the area pseudonorm
$\rho^+ + \rho^-$ is an  unbounded  norm 
  on $\Ham$.
\end{cor}
 \begin{proof}  See  Lemma~\ref{le:rho}.
\end{proof}

This result improves
\cite[Thm~1.2]{Mcv} which established the nontriviality of
$\rho^+ + \rho^-$ only for the cases 
$\om|_{\pi_2(M)} =0$ and $M=\C P^n$. 

The real problem in understanding the onesided pseudonorms
$\rho^{\pm}$ is caused by the possible existence of {\it short loops}, i.e. loops
in $\Ham$ that are generated by Hamiltonians for 
which $\rho^+$ (or $\rho^-$) is small; see \cite[\S1.4]{Mcv}.  These are still not understood in the closed case in general, but, as we now explain, it turns out that
they cause no problem when  $M$ is noncompact. 

Let $(M,\om)$ be a noncompact manifold without boundary, and 
 $U\subset M$ be open with compact closure.
  Denote by
$\Ham^c\, U$  the group of symplectomorphisms 
generated by functions $H_t$
with support in $U$ and by 
$\THam\,\!^c\, U$ its universal cover.   Denote by $\Trho_U^+$ the positive part of the Hofer norm on $\THam\,\!^c\,U$ and by $\rho_U^+$ the induced function on $\Ham^c\,U$.   Notice that in principle $\Trho_U^+$  might depend on $U$. But clearly $U'\subset U$ implies $\Trho_{U'}^+\ge \Trho_{U}^+$. A similar 
remark applies to $\rho_U^+$.

The following result was suggested by a remark  in an early version of
 \cite{EP6}. Note that $(M,\om)$ can be arbitrary here; in particular it need satisfy no special conditions at infinity. 

\begin{prop}\label{prop:rho}  Suppose that $(M,\om)$ is noncompact and that $U$ is an open subset of $M$ with compact closure.   
Then for all  $H:M\to \R$ with  support in $U$ there is 
$\de= \de(H,U)>0$ such that $\rho_U^+(\phi_t^H) = t\max H$ 
for all $0\le t\le \de$.
In particular
$\rho_U^++\rho_U^-$ is  a nondegenerate norm on $\Ham^c\, U$.
\end{prop}

The proof is given at the end of \S\ref{s:sp}. 
 Note that $\THam\,\!^c\,U$ always has infinite Hofer diameter because of the existence of the Calabi homomorphism. However, even though
 $\THam\,\!^c\,U$ is not a simple group,
  the kernel of the Calabi homomorphism
 is simple.  Hence the second statement in the above proposition follows from the first because the element  $\phi_t^H$
 belongs to this kernel when  $\int_U H\om^n=0$.
 For a brief discussion of 
other issues that arise in the noncompact case, 
see \cite[Remark~3.11]{Mcsn}.

Finally we describe another  class of  examples, 
the one point blow ups.   These have
 $N|(n-1)$  but may be chosen to satisfy the other conditions of 
Theorem~\ref{thm:main2}.

\begin{prop}\labell{prop:bl}  {\rm (i)}  Let $M$
 be a sufficiently small one point blow up of any closed symplectic manifold $(X,\om_X)$ such that at least one of  $[\om_X], c_1(X)$ 
does not vanish on  $\pi_2(X)$.
Then the asymptotic spectral invariants do not descend to $\Ham$.\SSS

\NI {\rm (ii)}   If $M$ is 
the one point blow up of a closed symplectic $4$-manifold $X$
such that both  $[\om_X]$ and $ c_1(X)$ 
 vanish on  $\pi_2(X)$
 then the asymptotic spectral invariants do descend to $\Ham$.
\end{prop} 

The proof is given
in \S \ref{ss:ex}.

\begin{cor}  If $M$ is as in Proposition~\ref{prop:bl}(ii) then
$\Ham$ supports  a nontrivial Calabi quasimorphism.
\end{cor}

The proof is contained in the following remark.

\begin{rmk}\labell{rmk:qm}\rm  Other potential applications of these results arise from the work of  Entov--Polterovich \cite{EP1,EP2,EP3,EP4}.  
 They denote the asymptotic spectral invariant given by an
 idempotent $e$ in $\Aa_M: = QH_{2n}(M)$ by
$$
\mu_e:\THam\to \R,\quad \mu_e(\cdot): = \oc(e,\cdot).
$$ 
It is immediate that $\mu_e$ descends to $\Ham$ iff
$\mu_e|_{\pi_1(\Ham)}$ vanishes; see Proposition~\ref{prop:des}.   
The Entov--Polterovich results about nondisplaceablity (see \cite{EP3} for example) apply whether or not
$\mu_e$ descends to $\Ham$. 
But if one is interested in questions about the structure of the Hamiltonian group itself, for example what quasimorphisms or norms it might have or what discrete subgroups it might contain, then our results are relevant.

If  $e$ is an idempotent such that
$e\Aa_M$ is a field,\footnote
{
When $\Aa_M$ is semisimple, it is tempting to think that $\mu_e$ is a quasimorphism for every idempotent, and in particular for $e=\1$.  However, this is not true, as is shown  by the example of the small one point blow up of $\C P^2$.  The calculations in \cite{Mcv} (see also \cite{Os2}) 
show that there is an element $Q\in \Aa_M$ such that $\nu(Q^k) + \nu(Q^{-k})\to \infty$. (See \S2 for notation.)  Further there is a constant $c$ and an element $\al\in \pi_1(\Ham)$ such that
$\mu_{\1}(\al^k): = \nu(Q^k) + kc$ for all $k\in \Z$.  But if $\mu_{\1}$ were a
quasimorphism it would restrict to a  homomorpism on the abelian subgroup  
$\pi_1(\Ham)$  and we would have
$\nu(Q^k) + \nu(Q^{-k})=0$.
  Hence it cannot be a
 quasimorphism.  Rather it is related to the maximum of two different quasimorphisms:
since $\1=e_1+e_2$ where each $e_i$ is minimal,   $\mu_{\1} \le \max(\mu_{e_1}, \mu_{e_2})$ with equality at all elements where 
the $\mu_{e_i}$ take different values; cf. equation (20) in \cite{EP3}.
}
 then
 Entov--Polterovich show that $\mu_e$ is a homogeneous
quasimorphism, i.e. there is $C>0$ such that for all $k\in\Z$ and
$\Tphi,\Tpsi\in \THam$
$$
\mu_e({\Tphi}\,\!^k) = k\mu_e(\Tphi),\qquad |\mu_e(\Tphi\,\Tpsi)-\mu_e(\Tphi)-\mu_e(\Tpsi)|\le C.
$$
In particular, the restriction of $\mu_e$ to the abelian subgroup
$\pi_1(\Ham)\subset \THam$ is  a homomorphism.  
Therefore,  $\mu_e$ descends to a quasimorphism on $\Ham$ exactly if
it vanishes on $\pi_1(\Ham)$. It has the Calabi property of 
 \cite{EP1} by construction.

There are rather few known manifolds
(besides $\C P^n$)  for which $\Aa_M$ contains an idempotent $e$
such that both
 $e\Aa_M$ is a field and $\mu_e$ 
descends.  For these conditions work in opposite directions; we need many nontrivial Gromov--Witten invariants for $e\Aa_M$ to be a field, but not too many (or at least  the Seidel representation of \S\ref{s:sdpr} should be controllable) if $\mu_e$ is to descend.
In this paper  
 the one set of new examples with a suitable idempotent $e$ are those of
Proposition~\ref{prop:bl} (ii).  Indeed, the calculations in
\cite[\S2]{Mcu} (see also \cite{EP4})  show that in this case there is an idempotent $e\in \Aa_M$ such that $e\Aa_M$ is a field, while we show here that $\mu_e$ descends.

However, even if $e\Aa_M$ is not a field, Entov--Polterovich
show in \cite[\S7]{EP2} that 
 $\mu_e$ interacts in an interesting  way with  
the geometry of $M$.  
For example, it
 provides a lower bound for the so-called fragmentation norm\footnote
 {
 If $U$ is an open subset of $M$,   $\|\Tphi\|_U$ is defined to be the minimal number $k$ such that $\Tphi$ can be written as a product of $k$
  symplectomorphisms  each conjugate to  an element in $\THam\,\!^c(U)$, the universal cover of the group of compactly supported Hamiltonian symplectomorphisms of $U$.}
  $\|\cdot\|_U$ on $\THam$ of the form  
 $$
 |\mu_e(\Tphi\Tpsi) - \mu_e(\Tphi) - \mu_e(\Tpsi)|\le K \min\{\|\Tphi\|_U,\|\Tpsi\|_U\},\quad\mbox{ for all } 
 \Tphi,\Tpsi\in \THam,
 $$
 where $K$ is a constant that depends only on the open set $U$ and $U$ is assumed displaceable, i.e. there is $\psi\in \Ham\, M$ such that  $U\cap \psi(U)=\emptyset$.  If $e\, \Aa_M$ is a field then 
 the quantity on LHS is bounded and $\mu_e$ is a quasimorphism as in (ii) above.
 On the other hand, if this quantity 
 is unbounded then the fragmentation norm $\|\cdot\|_U$ is also unbounded on $\THam$.  Theorems \ref{thm:main1} and 
 \ref{thm:main2} 
 allow one to transfer these results to $\Ham$ in many cases;
 cf. Burago--Ivanov--Polterovich \cite[Ex.~1.24]{BIP}.\SSS
  \end{rmk}

 \begin{rmk}\labell{rmk:spectv}\rm (i) (Properties of the spectral norm and its variants.)\,
 Consider the function $V$ 
 on $\THam$ given by
 $$
V(\Tphi):  =c(\1,\Tphi) + c(\1,\Tphi\,\!^{-1}),
 $$
  and the corresponding  function $v$ induced on $\Ham$:
 $$
v(\phi): = \inf \bigl\{V(\Tphi)\,|\, \Tphi\mbox{ lifts }\phi\bigr\}.
$$
 Schwarz and Oh  showed that $v$ is  a conjugation invariant norm  on $\Ham$, called the {\it spectral norm}.
As noticed by Entov--Polterovich ~\cite{EP1},
this norm is bounded when $M=\C P^n$ or, more generally, when
$QH_*(M)$ is a field with respect to suitable coefficients.\footnote
{
Cf. Albers \cite[Lemma~5.11]{Al}.  In this context it matters which 
coefficients are used for quantum homology; compare the approaches of Ostrover~\cite{Os2} and Entov--Polterovich~\cite{EP4}.}
For these hypotheses imply that  $c(\1,\cdot)$ is an
 (inhomogeneous) quasimorphism on $\THam$, so that
 $$
 V(\Tphi) = |c(\1,id)-c(\1,\Tphi) - c(\1,\Tphi\,\!^{-1})|\le {\it const}.
 $$
 Note that Ostrover's argument does not apply here since 
$ V(\Tpsi_s)$ remains 
bounded  on the path $\Tpsi_s,s\ge 0,$ 
defined in equation (\ref{eq:sp}).

 The supremum of $v(\phi)$ for $\phi\in \Ham$ is called  the {\it spectral capacity}  of $M$, cf.
    Albers  \cite[equation~(2.47)]{Al}.  Since $V(\Tphi)\ge c(\1,\Tphi)-c(pt,\Tphi)>0$ for $\Tphi\ne id$, 
    this capacity can be thought of as  a measure of spectral spread. 
 When $M$ is the standard torus $T^{2n}$ it is well 
 known that there are (normalized) functions $H:M\to\R$  whose flow $\Tphi\,\!^H_t$ has the property that $V(\Tphi\,\!^H_k) = k V(\Tphi\,\!^H_1)\ne 0$; see
the discussion  after Question~8.7 in  \cite{EP2}.
Since Theorem~\ref{thm:main2} implies that $V=v$ in this case, the spectral capacity of $T^{2n}$ is infinite.  In general  the spectral capacity is poorly understood. For example, it is not known whether it is always infinite when $QH_*(M)$ is very far from being a field,
for example if $(M,\om)$ is  aspherical.
  
As we explain in the proof of Lemma~\ref{le:rho}, 
$c(\1,\Tphi)\le \rho^+(\Tphi\,\!^{-1})$.  Hence one might think of
 minimizing
$c(\1,\Tphi)$ and $c(\1,\Tphi\,\!^{-1})$ separately over the lifts of $\phi$ as we did for the Hofer norm.  However, in general 
this procedure gives nothing interesting
 since $c(\1,\Tphi)$  certainly can be negative and may well have
no lower bound over a set of lifts.  Instead, one can 
consider the  functions
$$
\ovV: \THam\to \R,\quad \Tphi
\mapsto \oc(\1,\Tphi) + \oc(\1,\Tphi\,\!^{-1}),
$$
and
$$
\ovv: \Ham\to \R,\quad
\phi\mapsto \inf \,\bigl\{\ovV(\Tphi)\,|\, \Tphi\mbox{ lifts }\phi\bigr\}.
$$
(Thus $\ovV$ is a symmetrized version of the function $\mu_{\1}$ considered in the previous remark.)
If $c(\1,\cdot)$ is a quasimorphism then $\oc(\1,\cdot)$ is a 
homogeneous quasimorphism and hence satisfies
$\oc(\1,\Tphi)=
-\oc(\1,\Tphi\,\!^{-1})$.  Therefore, in this case, $\ovV\equiv 0$.
 On the other hand, as we pointed out in 
Remark~\ref{rmk:qm}, $\ovV$ does {\it not} vanish on 
$\pi_1(\Ham)\subset \THam$ when
$M$ is a small blow up of $\C P^2$.   Our remarks above 
 imply that $T^{2n}$  has 
infinite $\ovV$-diameter, but again very little is known about $\ovV$
for general $M$.
 
 Although these variants of $v$ have some uses,  they
are unlikely to be (pseudo)norms, since, as we explain in
Remark~\ref{rmk:asymp} (iii),
 they probably never have the property $m(fg)\le m(f) + m(g)$. One might also think of replacing the class  $\1$ by some other idempotent $e$. But it is easy to see that
  $\nu(e)>0$ for any such $e$.  Hence the resulting function would not take the value $0$ at the identity $\id\in \THam$. \MS

\NI (ii)  The proof of Corollary~\ref{cor:rho+} 
 compares $\rho^+ + \rho^-$ with the Schwarz--Oh norm $v$.
Since there are $\Tphi \ne \1$
such that $c(\1,\Tphi) \le 0$ (see equation (\ref{eq:sp}))  this approach does not help with the other  problem left open in \cite{Mcv}, namely the question of whether   $\rho^+$  ever vanishes on some $\phi\ne \1$.
\end{rmk}

\NI {\bf Acknowledgements}  Many thanks to Leonid Polterovich for very helpful comments on an earlier draft of this paper, and also to Peter Albers, Alvaro Pelayo and the referee for detailed comments that have helped improve 
the clarity of the exposition. 

\section{Spectral invariants}\labell{s:sp}

In this section we discuss the basic properties of spectral invariants 
and prove Lemmas~\ref{le:inf} and \ref{le:rho} concerning
the Hofer diameter of $\Ham$, as well as Proposition~\ref{prop:rho}.  We assume throughout that $M$ is closed unless explicit mention is made to the contrary.

\subsection{Quantum homology}\labell{ss:qh}

To fix notation we list some facts about  the small quantum homology $QH_*(M): = H_*(M)\otimes \La$.
We shall take  coefficients  $\La: = \La^{univ}[q,q^{-1}]$, where $q$ is a variable of degree $2$ and $\La^{univ}$ is the field\footnote
{
We use the ground field $\R$ here since later on we use homology with $\R$ coefficients, but could equally well take $r_i\in \C$
as in Entov--Polterovich~\cite{EP3}.}
 of generalized Laurent series in $t^{-1}$ with elements 
 $$
 \la = \sum_{i\ge 0} r_i t^{\eps_i}, \quad r_i,\, \eps_i\in \R, \
 \eps_i > \eps_{i+1},  \ \eps_i\to -\infty.
 $$
 Thus $\La^{univ}$ has a valuation $\nu:\La^{univ}\to [-\infty,\infty)$ given by
 $\nu(\la): = \max\{\eps_i: r_i\ne 0\}$.  Observe that $ \nu(0) = -\infty,$
 and for all $\la,\mu\in \La$
 $$
\nu(\la+\mu)\le\max\bigl(\nu(\la),\nu(\mu)\bigr),\quad
 \nu(\la\mu) = \nu(\la) + \nu(\mu),\quad \nu(\la^{-1}) = -\nu(\la).
 $$
 This valuation extends to
  $QH_*(M)$ in the obvious way: namely,  for any $a_i\in H_*(M)$ we set
  $$
  \nu\left(\sum a_i \otimes q^{d_i} t^{\eps_i}\right) =
  \max\{\eps_i: a_i\ne 0\}.
  $$
   The quantum product $a*b$ of the elements $a,b\in 
 H_*(M)\subset QH_*(M)$ is defined as follows.  
 Let $\xi_i, i=0,\dots,m,$ be a basis for $H_*(M)$ with dual basis $\{
 \xi_i^{*_M}\}$.  Thus
 $$
 \xi_i^{*_M}\cdot_M\xi_j = \de_{ij}.
 $$
 We use this slightly awkward notation to
 reserve   $\xi_i^*$ for later use; cf. equation (\ref{eq:basis}). Also $\cdot_M$ (which is often simplified to $\cdot$) denotes the
 intersection product 
 $$
 H_d(M)\otimes H_{2n-d}(M)\to H_0(M)\equiv \R.
 $$
 Further, denote by $H_2^S(M)$ the spherical homology group, i.e. the image of the Hurewicz map $\pi_2(M)\to H_2(M)$. 
 Then 
\begin{equation}\labell{eq:QH}
 a*b: = \sum_{i,\be \in H_2^S(M)} \bla a,b,\xi_i\bra^M_\be\,\, \xi_i^{*_M}\,\otimes q^{-c_1(\be)} t^{-\om(\be)},
\end{equation}
 where $\bla a,b,\xi_i\bra^M_\be$ denotes the Gromov--Witten invariant   that counts curves in $M$ of class $\be$  through the homological
 constraints
 $a,b,\xi_i$.
 Note that 
 $\deg(a*b) = \deg a + \deg b - 2n$, and the identity element is $\1: = [M]$.
 The product is extended to $H_*(M)\otimes \La$ by linearity over $\La$.

Later it will be useful to consider the  $\La$-submodule
$$
\Qq_-: = \bigoplus_{i<2n} H_i(M)\otimes \La.
$$  
(Here as usual $2n = \dim M$.)   The following (easy) result was 
proved in
\cite{Mcu}.

\begin{lemma}\labell{le:Qq} {\rm (i)} 
If $M$ is not strongly uniruled then $\Qq_-$ is an ideal in $QH_*(M)$.
 \SSS
 
 \NI {\rm (ii)} If $\Qq_-$ is an ideal and if $u\in 
 QH_{2n}(M)$ is invertible, then 
 $u = \1\otimes \la + x$ for some nonzero $\la\in \La$ and some
 $x\in \Qq_-$.
 \end{lemma}
 
Note that if   quantum multiplication is undeformed  then the elements of
$\Qq_-$ are nilpotent so that  all  elements 
$\1\otimes \la + x$ with $\la\ne 0$ and $x\in \Qq_-$ are invertible.
However,  if we assume only that
 $\Qq_-$ is an ideal, then an element of the form
$\1\otimes \la + x$   might not be invertible. For example  it could  be a nontrivial idempotent.
 See \cite{Mcu} for further details.

 \NI{\it
 Note:} whenever we write a unit as $\1\otimes \la + x$ we assume, unless explicit mention is made to  the contrary, that 
 $\la\ne 0$ and $x\in \Qq_-$.
 
\subsection{Spectral invariants and norms}\labell{s:si}
 
One way to estimate the length of a Hamiltonian path is to use the Schwarz--Oh {\it spectral invariants}; see \cite{Sch,Oh} and Usher~\cite{U}. (They are also explained in \cite[Ch~12.4]{MS2}.)
For each element $\Tphi\in \THam$ and each nonzero 
element $a\in QH_*(M)$
the number $c(a,\Tphi)\in \R$ has the following properties:
\begin{eqnarray}\labell{eq:1}
&&-\|\Tphi\|\; \le\; c(a,\Tphi) \;= \; c(a,\Tpsi\,\Tphi\,{\Tpsi}\,\!^{-1})\;\le\; \|\Tphi\|\;\;\mbox{ for } a\in H_*(M), \Tpsi\in \THam 
\\\labell{eq:2}
&& c(\la a,\Tphi) \;= \; c(a,\Tphi) +\nu(\la)\;\;\mbox{ for all }\la\in \La,\\\labell{eq:3}
&& 
c(a,\Tphi\circ\ga)\;=\; c(\Ss(\ga)*a,\Tphi)\;\;\mbox{ for all } \ga\in \pi_1(\Ham),\\
\labell{eq:4}
&& 
c(a*b,\Tphi\circ\Tpsi)\;\le\; c(a,\Tphi) + c(b,\Tpsi)
\;\;\mbox{ for all } a,b\in QH_*(M),  \Tphi, \Tpsi\in \THam.
\end{eqnarray}
The third property explains how these numbers
 depend on the path $\Tphi$.  Here, the element 
$\Ss(\ga)\in QH_*(M)$ is called  the {\it Seidel element of the loop} $\ga$
(see ~\cite{Sei,Mcq}). 
It is an invertible element of degree $2n=\dim M$ in $QH_*(M)$; we give a brief 
 definition in \S\ref{s:sd} below.  Further,  (\ref{eq:1})  implies that 
 $c(a, id) = 0$ for all $a\in H_*(M)$, where $id$ denotes the constant loop
 at the identity.  Hence, for all $\ga\in \pi_1(\Ham)$ 
 $$
 c(\1,\ga) = c(\Ss(\ga),id) = \nu(\Ss(\ga)),\qquad
 \oc(\1,\ga) = \lim_{k\to\infty} \frac{\nu(\Ss(\ga^k))}{k}.
 $$

\begin{rmk}\labell{rmk:asymp}\rm (i)
Equation 
(\ref{eq:1}) implies
that the lim inf defining the asymptotic invariants 
$\oc(a,\cdot)$ always exists. When $a^2=a$ standard arguments based on
(\ref{eq:4})
show that one can replace the lim inf by an ordinary limit;
 cf. \cite[\S4.2]{EP1}.  \MS

\NI (ii) The first two properties (\ref{eq:1}), (\ref{eq:2})
of the spectral invariants are obviously inherited by the asymptotic invariants.  The third is too, because $\pi_1(\Ham)$ lies in the center of $\THam$; cf. the proof of Proposition \ref{prop:des} below.
Moreover,  if 
$\Tphi, \Tpsi$ commute then the fourth property is also inherited by $\oc$. 
(This is the basis for the discussion in \cite{EP2} of {\it partial symplectic quasi-states}.) \MS

\NI (iii)
If $e$ is an idempotent it is easy to see that
$
\oc(e,\Tphi\Tpsi) \le c(e,\Tphi)+c(e,\Tpsi)$ for all $\Tphi,\Tpsi$.  But in general 
the asymptotic invariants only have good algebraic properties when
$e$ is an idempotent such that $e\Aa_M$ is a field; cf. Remark~\ref{rmk:qm}. However even in this case it is impossible to have
$\oc(e,\Tphi\Tpsi) \le \oc(e,\Tphi)+\oc(e,\Tpsi)$  for all $\Tphi, \Tpsi$. Indeed, as Polterovich points out,\footnote
{Private communication.} because  $\oc(e,\Tphi\,\!^{-1}) = -\oc(e,\Tphi)$ for such $e$,  if this inequality did always  hold one would have
$$
- \oc(e,\Tphi\Tpsi)\; = \;\oc(e,(\Tphi\Tpsi)^{-1})\; \le\; \oc(e,\Tphi\,\!^{-1}) +\oc(e,\Tpsi\,\!^{-1})\; =\; - \oc(e,\Tphi)-\oc(e,\Tpsi).
$$
But then we would have equality, i.e.
 $\oc(e,\cdot)$ would be a surjective homomorphism. 
 It would then descend to a
 nontrivial homomorphism $\Ham\to \R/\Ga$, where $\Ga$ is the image 
 of the countable\footnote
 {
 A proof of this classical result is sketched in \cite[Remark~9.5.6]{MS2}.} group  $\pi_1(\Ham)$.  Since $\Ham$ is perfect,
this  is impossible. 
(Incidentally, this argument also shows  that when $e\Aa_M$ is a field it is impossible that $c(e,\Tphi) = \oc(e,\Tphi)$ 
for all $\Tphi$.
This is obvious from the explicit calculation in formula (\ref{eq:sp}) below. On the other hand, the current argument is structural and hence might apply in other situations as well.)
\end{rmk}

Much of  following proposition
 is implicit in Entov--Polterovich \cite{EP1}.\footnote
 {
 In the published version, 
 it is claimed that (i:a) below is equivalent to the statement that just one invariant $c(a,\cdot)$ descends, and  the same claim is made for (ii:a).   However, this could be false, though I am not aware of any explicit counterexample at this time.
 Corrected versions of (i:b) and (ii:b) appear below; they are not used elsewhere in the paper.}

\begin{prop}\labell{prop:des} {\rm  (i)}  The following conditions are equivalent:
\SSS

{\rm (a)}
All spectral invariants $c(a,\cdot),a\in QH_*(M),a\ne 0,$  descend to $\Ham$;\SSS

{\rm (b)} 
The spectral invariant $c(\1,\cdot)$  descends to $\Ham$;\SSS

{\rm (c)}
$c(\1,\ga) = \nu(\Ss(\ga)) = 0$ for all $\ga\in \pi_1(\Ham)$;\SSS

{\rm (d)}  For all $\ga\in \pi_1(\Ham)$, 
$\Ss(\ga) =  \1\otimes \la +x$ where  $\nu(\la) = 0$ and  $\nu(x)\le 0$.\MS

\NI
{\rm  (ii)} The following conditions are equivalent:
\SSS

 {\rm (a)}   All asymptotic spectral invariants 
$\oc(a,\cdot),a\in QH_*(M), a\ne0,$ descend to $\Ham$;
\SSS

{\rm (b)}
The  asymptotic spectral invariant
$\oc(\1,\cdot)$ descends to $\Ham$;\SSS

{\rm (c)}  For all $\ga\in \pi_1(\Ham)$, 
$$
c(\1,\ga)\ge0,\;\mbox{ and }\; \oc(\1,\ga) = \lim_{k\to\infty}\;\frac{\nu(\Ss(\ga^k))}k\ = \;0.
$$
\NI {\rm (iii)}  If the spectral invariants descend then so do the asymptotic invariants.
\end{prop}
\begin{proof}    Consider (i).  The equivalence of the first 
three conditions is immediate from (\ref{eq:2}), (\ref{eq:3}) and
(\ref{eq:4}).
Part (ii) of the next lemma shows that (c) implies (d).  Conversely,  suppose that (d) holds, i.e. all the $\Ss(\ga)$ may be written as $\1\otimes \la + x$ where
 $\nu(\la) = 0\ge\nu(x)$.  Then $\nu(\Ss(\ga))\le 0$ for all $\ga$.
 Moreover, if $\nu(\Ss(\ga))< 0$ for some $\ga$ then 
 $$
0 \ = \ \nu(\1) \ = \ \nu\bigl(\Ss(\ga)\,\Ss(\ga^{-1})\bigr) \ \le 
\  \nu\bigl(\Ss(\ga)\bigr)+ \nu\bigl(\Ss(\ga^{-1})\bigr) \ < \  0
 $$ 
which is a contradiction.  Hence (c) holds.

 Now consider (ii).   
If the invariants $\oc$ descend then in particular we must have
$\oc(\1,\ga) = \oc(\1,id) = 0$ for all $\ga$.   Further, if
$c(\1,\ga)= -\eps < 0$ then $c(\1, \ga^k)\le -k\eps$ which implies that
$\oc(\1,\ga)< 0$. 
Therefore (a) implies (b) implies (c).
  To prove that (c) implies (a), observe that
because $\ga\in \pi_1(\Ham)$ is in the center of $\THam$,
 the subadditivity relation (\ref{eq:4})  for $c$ implies that
$$
\oc(a, \Tphi\circ\ga) \;= \;
\liminf {\textstyle \frac1k}\bigl(c(a,\Tphi\,\!^k\ga^k)\bigr) \;\le\; \oc(a,\Tphi) + \oc(\1,\ga).
$$
Hence if $\oc(\1,\ga)=0$ always,
$$
\oc(a, \Tphi\circ\ga)\; \le\; \oc(a, \Tphi)\;= \;\oc(a, (\Tphi\circ\ga)\circ\ga^{-1})\;\le\; \oc(a, \Tphi\circ\ga).
$$
Thus we must always have $\oc(a, \Tphi\circ\ga)=\oc(a, \Tphi)$.
Therefore (a), (b) and (c) are equivalent. 
This proves  (ii);  (iii) is immediate.
\end{proof}

\begin{lemma}\labell{le:help} {\rm (i)} Every element of the form $u=\1\otimes \la +y$ with  $\nu(\la) \ge\nu(y)$ is invertible in $QH_*(M)$.
\SSS

\NI {\rm (ii)} Let $u\in QH_*(M)$ be an invertible element of the form
 $u=\1\otimes \mu + y$ where  $\mu$ might be zero, and suppose that
 $\nu(u) = 0$.  Then either $\nu(\mu) = 0 = \nu(u^{-1})$ or $\nu(u^{-1})> 0$.
 \end{lemma}
\begin{proof}
First note that by the definition of the quantum product 
\begin{equation}\label{eq:*}
  \nu(x*z-x\cap z) \le \nu(x)+\nu(z)-\de \quad \mbox{ for all } x,z\in QH_*(M),
\end{equation}
where $\de> 0$ is the minimal energy $\om(\be)$ of a class $\be\ne0$
with nonzero Gromov--Witten invariant $\bla a_1,a_2,a_3\bra_\be$.
Hence, because all elements in the undeformed ring  $H_{<2n}(M)$ are nilpotent of order $\le n+1$,  
$$
\nu(a^{n+1}) \le -\de \qquad \mbox{for all } a\in H_{<2n}(M).
$$
(Here and subsequently $a^k$ denotes the $k$-fold quantum product
$a*\dots*a$.)
Hence if $z\in \Qq_-$ has $\nu(z)\le 0$ then
$
\nu(z^{(n+1)k}) \le -k\de\to -\infty
$
as $k\to \infty$.  Because   $\nu(z^i)$ is a nonincreasing sequence when $\nu(z)\le 0$, it also diverges to $-\infty$.   Thus 
$\sum_{i\ge 0} z^i$ is a well defined element of $QH_*(M)$.
Hence
  $\1\otimes \la+y$ is invertible with inverse 
$$
u^{-1}: = \
\la^{-1}(\1+\la^{-1}y)^{-1} \ = \ \la\sum (-1)^i (\la^{-1}y)^i \ \in  \ QH_*(M).
$$
This proves (i).

To prove (ii)  suppose that $u=\1\otimes \mu + y$ is a unit  with $\nu(u) =0$.  Because $\1$ and  the components of $y$ are linearly independent, $\nu(u) = \max\bigr(\nu(\mu), \nu(y)\bigr)$.  Similarly, if we write
$u^{-1} = \1\otimes \la + x$, we have $\nu(u^{-1}) = 
\max\bigr(\nu(\la), \nu(x)\bigr)$.  

If $\nu(\mu) = 0$ then the above formula for $u^{-1}$ shows that $\nu(u^{-1}) = 0.$  So suppose that
$\nu(\mu)<0$.
Because  $u\,u^{-1} = \1$, at least one of the terms $1\otimes \la\mu$ and $xy$ must contain $\1$ with a nonzero rational coefficient.   If it is the former then $\nu(\la)=-\nu(\mu)>0$, while if it is the latter
then we must have $\nu(x) + \nu(y)\ge \de$ so that $\nu(x)\ge \de$.
In either case $\nu(u^{-1})>0$.
 \end{proof}

\begin{cor}\labell{cor:des}
The asymptotic spectral invariants descend to $\Ham$ if
for all $\ga$ there is $m$ such  that 
$\Ss(\ga^m)= \1\otimes \la +x$ where $\nu(\la) = 0$ and $x$ is nilpotent.
\end{cor}
\begin{proof} Suppose that $x^N=0$.
Let $K = \max \{\nu(x^k)\,|\, 0\le k\le N\}.$
Then because $\Ss$ is a homomorphism and $\nu$ is subadditive
$$
\nu\bigl(\Ss(\ga^{mk})\bigr)= \nu\big((\Ss(\ga^m))^{k}\bigr)
= \nu(\la^k + k\la x + \dots) \;\le\; \max\{
k\nu( \la), \nu(\la x), \dots\}\;\le\; K,
$$
  for all $k\in \Z$ and all $\ga$. 
Hence $\oc(\1,\ga^m) = 0$, and hence $\oc(\1,\ga) = 0$.
\end{proof}

\begin{rmk}\rm
It is not hard to find conditions  under which
every invertible element in $QH_*(M)$ has the form $\1\otimes \la +x$
where $x$ is nilpotent.  For example, as we noted above, this is always true if the quantum multiplication is undeformed.
For other cases see Lemma~\ref{le:cN}. 
Therefore the main difficulty in 
showing that the (asymptotic) spectral invariants 
descend lies in ensuring that the condition $\nu(\la)=0$ holds. 
\end{rmk}

The numbers $c(a,\Tphi)$ are defined by looking at the filtered Floer complex of the generating 
mean normalized
 Hamiltonian $H$, and  turn out to be particular critical values\footnote
{
In fact this property is essential to the existence of the spectral invariants as functions on $\THam$; the spectral invariants  
are first defined as functions on the space of Hamiltonians $H$
and one needs the spectrality property to conclude that they 
actually depend only on the element in $\THam$ defined by the flow of $H$.
}
 of the corresponding
action functional $\Aa_H$; cf. Oh \cite{Oh1} and Usher \cite{U}.
  Thus each  $c(a,\Tphi)$ corresponds to a particular fixed point of the endpoint $\phi_1\in \Ham$ of $\Tphi$.  

It is usually very hard to calculate them. 
However, if
 $\Tphi: = \Tphi^H$ is generated by a $C^2$-small  mean normalized  Morse function $H:M\to \R$, then
\begin{equation}\labell{eq:CM}
  c(a,\Tphi^H)= c_M(a,-H),\qquad a\in QH_*(M)
\end{equation}
where $c_M(a,-H)$ are  the corresponding invariants obtained from the Morse complex of $-H$.  
 These  are defined as follows.  Denote by
 $CM_*(M,K)$
  the usual Morse complex for the Morse function $K$.
  For each $\ka\in \R$, it has a subcomplex
 $CM_*^\ka(M,K)$   generated by the
  critical points $p\in \Crit(K)$ with critical values $K(p)\le \ka$.
 Denote by $\io_\ka$ the inclusion of the homology  $H_*^\ka$
 of this subcomplex  into $H_*^{\infty} \cong H_*(M)$. Then 
 for each $a\in H_*(M)$ 
\begin{equation}\labell{eq:M}
 c_M(a,K): = \inf\{\ka: a\in {\rm Im}\,\io_\ka\}.
\end{equation}

\NI {\bf Ostrover's construction.}
As pointed out by Ostrover~\cite{Os}, one can use the continuity properties of the $c(a,\Tphi)$ to find a path $s\mapsto \Tpsi_s$ in $\THam$ whose spectral invariants tend to $\infty$ as $s\to \infty$.

Normalize $\om$ so that $\int_M\om^n=1$.  Let $H$ be a small mean normalized  Morse function, choose an 
open set $U$ 
that is displaced by $\phi_1^H$ (i.e. $\phi_1^H (U)\cap U=\emptyset$)
and let $F:M\to \R$ be a function with support in $U$ and with nonzero integral $I: = \int_M F\om^n$, so that $F-I$ is mean normalized.
 Denote the flow of $F$ by $f_t$ and consider 
the path $\Tpsi_s: = \{f_{ts}\phi_t^H\}_{t\in [0,1]}$ in $\THam$ as $s\to \infty$.  For each $s$, $\Tpsi_s$ is generated by the Hamiltonian
$$
F_{st}\# H: = sF + H\circ f_{st}.
$$
The corresponding mean normalized Hamiltonian
is
$$
K_s: = sF+ H\circ f_{st} - sI.
$$ 
By construction, $f_{s}\phi_1^H$ has the same fixed points
as $\phi_1^H$, namely the critical points of $H$.  Hence the continuity and spectrality  properties of $c(a,\Tpsi_s)$ imply that for each $a$ the fixed point $p_a$ whose critical value is
$c(a,\Tpsi_s)$  remains {\it unchanged} as $s$ increases.
But the spectral value does change. In fact, if $a\in H_*(M)$, then when $s=0$  there is a critical point $p_a$ of $H$ such that
$c(a,\Tpsi_0) = c_M(a,-H) = -H(p_a).$  Hence 
\begin{equation}\label{eq:sp}
c(a,\Tpsi_s) = -K_s(p_a) =  -H(p_a) + sI,\;\;\mbox{ for all } s\in \R, a\in H_*(M).
\end{equation}
By ({\ref{eq:1}) it follows that $\THam $ has infinite Hofer diameter. 

The next result is well known.

\begin{lemma}\labell{le:inf} $\Ham$ has infinite diameter if either a 
spectral number or an asymptotic spectral  
 number descends to $\Ham$.
 \end{lemma}
 \begin{proof}  Let $\pi:\THam\to \Ham$ denote the projection.
     We will show that, for the above path $\Tpsi_s$, 
     $$
     \|\pi(\Tpsi_s)\| \;= \; \|\psi_s\|\;\to\;\infty \mbox{ \; as\; }s\to \infty.
     $$
       By definition, for each $s$ 
 $$
 \|\psi_s\| \;= \; \inf  \{\|\Tpsi_s\circ\ga\|: \ga\in \pi_1(\Ham)\}\;\;\ge\;\;
 \inf \{c(a,\Tpsi_s\circ\ga): \ga\in \pi_1(\Ham)\}
 $$
 where the inequality follows from equation (\ref{eq:1}).  If the spectral number $c(a,\cdot)$ descends to $\Ham$ then $c(a,\Tpsi_s\circ\ga)=c(a,\Tpsi_s)$
 for all $\ga$ so that the result follows from equation (\ref{eq:sp}).

Now suppose that the asymptotic spectral number $\oc(e,\cdot)$ 
 descends to $\Ham$ for some idempotent $e$.  We first claim that
 there are elements $\Tg_i, i=1,\dots,k-1,$ in $\THam$ that are 
 conjugate to $\Tphi^H_1$ and such that
   $\Tpsi_{sk}\Tg_1\cdots \Tg_{k-1} = (\Tpsi_s)^k $.
To see this, denote $a: = \Tilde f_s, b: = \Tphi^H_1$ so that
$\Tpsi_{sk} = a^kb$ and $ (\Tpsi_{s})^k = (ab)^k$. Then use the fact that
$$
(ab)^k = a^k b\cdot b^{-1} a^{-k} ab\dots ab =
a^k b\cdot b^{-1}\Bigl((a^{-k+1}ba^{k-1}) (a^{-k+2}ba^{k-2})\dots (a^{-1}ba) \Bigr)b.
$$
Next, observe that by (\ref{eq:4}) we have
 $c(e,(\Tpsi)^k) \le k c(e,\Tpsi)$ for all $k>1$.
   Hence
   $$
   \oc(e,\Tpsi_s\circ\ga)\;= \;\lim {\textstyle \frac 1k }
   c(e,(\Tpsi_s\circ\ga)^k)\;\le\;
   c(e,\Tpsi_s\circ\ga)\;\le\; \|\Tpsi_s\circ\ga\|.
   $$
   On the other hand,
\begin{eqnarray*}
     \oc(e,\Tpsi_s\circ\ga)\; \;= \;\;  \oc(e,\Tpsi_s)
   &=& \lim {\textstyle \frac 1k} c(e,(\Tpsi_{sk}) \Tg_{k-1}^{-1}\cdots \Tg_1^{-1})\\
   &\ge & \lim {\textstyle \frac 1k} \Bigl(c(e,\Tpsi_{sk}) - \sum_{i=1}^{k-1}c(e,\Tg_i)\Bigr),\\
&\ge &  s I - \|\Tphi^H_1\|,
\end{eqnarray*}
 where the first  inequality uses the identity $c(e,fg) \ge c(e,f)-c(e,g^{-1})$
 which follows from (\ref{eq:4}),
and the second uses (\ref{eq:sp}) and the fact that the $\Tg_i$ are conjugate to $\Tphi_1^H$.   
 \end{proof} 

\begin{lemma}\labell{le:rho} The function  $\rho^+ + \rho^-$ is unbounded on $\Ham$ whenever the 
 asymptotic spectral numbers   descend to $\Ham$.
\end{lemma}

\begin{proof}  The proof of the inequalities in
(\ref{eq:1}) actually shows that 
$$
\int_0^1 \min (-H_t)\, dt\le c(\1,\Tphi) \le \int_0^1\max (-H_t)\,dt
$$
  for every mean normalized $H_t$ that generates $\Tphi$;
see for example \cite[Thm.12.4.4]{MS2}.
Now suppose 
 that $\oc(\1,\cdot)$ descends.  
Equation (\ref{eq:4}) implies that
$
\oc(\1,\Tphi\,\!^k)\le kc(\1,\Tphi)$.  Hence 
if $\Tphi$ is any lift of $\phi$,
 $$
 \oc(\1,\phi) = \oc(\1,\Tphi)  \le \int_0^1\max (-H_t)\,dt.
 $$
  Therefore, because $-H_{1-t}$ generates $\Tphi\,\!^{-1}$, we have
  $ \oc(\1,\phi) \le \rho^+(\phi\,\!^{-1})$.
 Applying this to the image $\psi_s\in \Ham$ of the element $\Tpsi_s$ of (\ref{eq:sp}) we find that
 $$
 \oc(\1,\psi_s) = s \le \rho^+(\psi_s\,\!^{-1}).
 $$
Thus, because $\rho^-\ge 0$, the function $\rho^+ + \rho^-$ is not identically zero and hence is a norm.  Moreover it is clearly unbounded.
 \end{proof}

\begin{rmk}\labell{rmk:cN}\rm 
 As remarked in the introduction, the Hofer pseudonorm 
 $\|\cdot\|$ is a norm on $\THam$ if and only if it restricts to a norm on $\pi_1(\Ham)$. The only way that I know to estimate $\|\ga\|$ for $\ga\in \pi_1(\Ham)$ is via the spectral invariant $c(\1,\ga) =
\nu(\Ss(\ga))$ which is a lower bound for $\|\ga\|$  by (\ref{eq:1}).
 If the spectral invariants descend to $\Ham$ then $c(\1,\ga)=0$ for all $\ga$ and one gets no information.   For example, there are no current methods to detect the Hofer lengths of the elements of
  $\pi_1(\Ham)$ when $(M,\om)$ is the standard $2n$-torus. (Of course, in this case  if $n>1$  it is also unknown  whether 
 $\pi_1(\Ham)$ itself is nonzero.)  On the other hand, there are manifolds such as 
 certain blow ups of $\C P^2$ and $S^2$-bundles over $S^2$ for which 
 $\pi_1(\Ham)$ is known and $\nu(\Ss)$ does restrict to an injective
 homomorphism; see \cite{Mcsn,Mcv} and the references cited therein. 
 \end{rmk}

We end this section by proving Proposition~\ref{prop:rho} about
the behavior of $\rho^+$ in the noncompact case.
We need to prove:

\begin{lemma}  For each $H:M\to \R$ with support in $U\ne M$ there is $\de=\de(H,U)>0$  such that  
$$
\rho^+_U(\phi_{tH}) = t\,\max H, \quad 0\le t\le \de.
$$ 
\end{lemma}
\begin{proof}  We will prove this for open sets $U$ with smooth boundary.
Since any precompact open set $U'$ is contained in such $U$  the result for $U'$ follows easily, since $\rho^+_{U'}$ is defined by taking an infimum over a smaller set than $\rho^+_U$.

Put a collar neighborhood $Y\times [-1,1]\subset M$ round the boundary of $U$
so that $$
  U\cap
 (Y\times [-1,1]) = Y\times [-1,0).
 $$
   Let $U_\la: = U\cup (Y\times [-1,\la])$.
    Choose any $\om$-compatible almost complex structure $J_0$ on  $M$ and choose $r>0$ so that for each $y\in Y\times \{1/2\}$ there is a symplectic embedding $f_y$ of the standard ball of radius $r$ into
    $Y\times [1/4,3/4]$ with center at  $f_y(0) =y$.
Then by monotonicity  there is $\de_0>0$ so that 
every $J_0$-holomorphic sphere
$u: S^2\to U_1$ whose image meets both $Y\times \{1/4\}$ and 
$Y\times \{1/2\}$  has energy $\ge \de_0$.

 One condition on the constant $\de$ is that $\de\|H\|< \de_0$, where $\|H\| = \max H - \min H$.  The other is that when $0<t\le \de$
 the function
  $t H$ should be sufficiently small in the $C^2$ norm for it to be possible to embed a
 ball of capacity $\|tH\|$ in $R_{tH}^-(\eps)$, the region \lq\lq under the graph" of $tH$.  This region is described in \cite[\S2.1]{Mcv}, and the ball embedding is constructed in \cite[II,~Lemma~3.2]{LM}.
 
Now let us suppose that  the lemma does not hold 
with this choice of $\de$.  Thus if $F: = t H$  for some $t\in[0, \de]$ we suppose that
 there is another Hamiltonian path $\phi^K_t$ in $\Ham^cU$ with $\phi^K_1=\phi^F_1$ that is
 generated by a Hamiltonian $K_t$  with support in $U$ and such that
 $$
 \int_0^1\int(\max_{x\in M} K_t) \;\om^n\;dt<   
\int(\max_{x\in M} F) \;\om^n.
 $$
 Under these conditions we show in \cite[\S2]{Mcv} that there is a Hamiltonian bundle 
 $$
 (U_1,\om) \to (P_1,\Om)\to S^2
 $$
  with the following properties:
 \MS
 
 \NI
{\rm (i)}  $P_1$ is trivial outside $U$.  More precisely,
    $(P_1,\Om)\to S^2$ contains a subbundle $(P,\Om)\to S^2$ with fiber $(U,\om)$ such that $(P_1\less P,\Om)\to S^2$ may be identified with the product bundle $((U_1\less U)\times S^2, \om\times \al)$ where $\al$ is   an area form on $S^2$.\SSS
    
     \NI
{\rm (ii)} $\int_{S^2}\al < \|F\|$.  \SSS
 
     \NI
{\rm (iii)}  $(P,\Om)$ contains a symplectically embedded ball $B$ of capacity $\|F\|$.
 \MS
 
 $P_1$ is constructed as the union $R_{K,F}(2\eps)$ of two bundles over $D^2$, one a smoothing of the region in $M\times [0,1]\times \R$  between the graphs of $K_t$ and of the  function $t\mapsto\max K_t$ and with anticlockwise  boundary monodromy $\{\phi^K_t\}_{t\in [0,1]}$, and the other a smoothing of the region between the graphs of $F$ and of the constant function $t\mapsto\min F$ and with 
 clockwise boundary monodromy 
 $\{\phi^F_t\}_{t\in [0,1]}$.  The bundle is trivial outside $U$ because $F,K$ have support in
  $U$.  Property (ii) is an immediate consequence of the construction.  The region below the graph of $F$ contains the ball, which is embedded near the maximum of $F$.

Now consider the Gromov--Witten invariant that counts 
holomorphic spheres in $P_1$ in the class  
$\si = [p\times S^2]\in H_2(P_1)$ for $p\in U_1\less U$ and through 
one point $p_0$.  If we restrict the class of allowable almost complex structures to those that are $\Om$-tame and equal to $J_0\times j$ outside 
$P$ and choose $p_0\in U_{1/4}$, then this invariant  is well defined  because the energy $\Om(\si)=\int_{S^2}\al$ of each curve is less than $\de_0$, the amount of energy needed for a curve to
enter the boundary region  $P_1\less P_{1/2} \cong (Y\times (1/2,1])\times S^2$. (Note that because $(P_1\less P, J_0\times j)$ is a product, the energy of any curve that enters this region is at least as much as the energy of its projection to a fiber.) 
We shall call this invariant
$\bla pt\bra^{P_1}_\si$, though in principle it might depend on
the product structure  imposed near the boundary of $P_1$.

Since there is a unique $\si$-sphere though each point in 
$P_{1/4}\less P$,
$\bla pt\bra^{P_1}_\si=1$. But then, by a standard argument,
the nonsqueezing theorem of \cite[\S1.4]{Mcv} holds for this fibration, i.e. every ball  of radius $r$ embedded in ${\rm Int\,} P$ has capacity $\pi r^2\le \Om(\si)$.
But this contradicts properties  (ii) and (iii) of $P_1$.

We conclude with a sketch proof of the nonsqueezing theorem.  Let $f:B_r\to {\rm Int\,} P_1$ be a symplectic embedding of the standard $r$-ball  with center $f(0) $ and image $B$. 
Choose an $\Om$-tame almost complex structure $J$
on $P_1$ that is normalized on $P_1\less P$ and equal to $f_*(j_n)$  on $B$, where $j_n$ is the standard structure on $\C ^n$.  
By hypothesis there is   a $J$-holomorphic sphere $C$ in class $\si$ through  $f(0)$, and so
$$
\Om(\si) = \int_ C\Om\; >\; \int_{C\cap B} \Om \;\ge\; \int_{f^{-1}(C\cap B)}\om_0 \;\ge\;  \pi r^2,
$$
where the last inequality holds because $f^{-1}(C\cap B)$ is a properly embedded holomorphic curve through the center of $B_r\subset \C^n$.
For more details, see \cite[Ch.~9]{MS2}.  
\end{proof}

\begin{rmk}\rm  Because the homology of $M$ has  no fundamental class   when $M$ is noncompact, it is not clear how to 
understand quantum homology and the Seidel representation in this case.
Nevertheless, as we will see below, for the problems under consideration
here we do not need to know everything about
the Seidel element $\Ss(\ga)$, but just some facts about 
the coefficient of the fundamental class $\1$ in $\Ss(\ga)$. This coefficient is given by 
  counting the number of sections of the corresponding fibration $P_\ga\to S^2$ that go through a point, i.e. by a Gromov--Witten invariant of the form
  $\bla pt\bra^{P_\ga}_\si$.  Therefore the above argument fits naturally into the framework developed below for the closed case.
\end{rmk}

\section{The main argument}\labell{s:sd}

This section explains the main ideas and proves 
Theorems~\ref{thm:main1} and \ref{thm:main2} modulo some calculations of Gromov--Witten 
invariants that are carried out 
in \S\ref{s:GW}.

\subsection{The Seidel representation}\labell{s:sdpr}

Our aim in this section is to prove the following result.

\begin{prop}\labell{prop:sch} {\rm (i)} If 
$\om$ vanishes on 
$\pi_2(M)$ then the spectral invariants descend to $\Ham$. \SSS

\NI {\rm (ii)}   If $\rank\, H_2(M) = 1$, $n+1\le N<\infty$ and
$(M,\om)$ is not strongly uniruled then  
the spectral invariants descend.  The same conclusion holds 
if $N\ge n-1$ and $(M,\om)$ is negatively monotone.\SSS

\NI {\rm (iii)}   If $\rank\, H_2(M) = 1$ and $n+1\le N<\infty$ then  
the asymptotic spectral invariants descend. 
\end{prop}

The main tool is the Seidel representation $\Ss$.  We will see that the stringent hypotheses above guarantee that the 
relevant moduli space of spheres 
has no bubbling so that it  is compact.  The argument goes back to 
Seidel\footnote{Private communication.} though it was first published by Schwarz,
 who showed in \cite[\S4.3]{Sch} that the spectral invariants descend when
both  $\om$ and $c_1$ vanish on 
$\pi_2(M)$.  We generalize this result in (i) above.  
Part (iii)  is a mild generalization of a result of  
Entov--Polterovich \cite{EP1} who proved that the
asymptotic spectral invariants descend when $M=\C P^n$.
 
After defining  $\Ss$, we shall
calculate it   under the conditions of the above proposition
(see Lemma~\ref{le:cN}).  Finally we prove the various cases of the proposition.

We shall think of the Seidel representation as
 a homomorphism 
 $$
 \Ss: \pi_1(\Ham(M,\om))\to QH_{2n}(M)^\times,
 $$
to the degree $2n$ multiplicative units  
of the small quantum homology ring, where  $2n: = \dim M$.
To define it, observe that  each loop $\ga=\{\phi_t\}$ in  $\Ham\,M$ gives rise to an 
$M$-bundle
$P_\ga\to S^2$ defined by the clutching function $\ga$:
$$
P_\ga: = (M\times D_+)\,\cup \,(M\times D_-)/\!\sim, \qquad (\phi_t(x),e^{2\pi it})_+\sim
(x,e^{2\pi it})_-. 
$$
(Here $D_\pm$ are two copies of the unit disc with union  $S^2$.) Because the loop $\ga$ is Hamiltonian, the fiberwise symplectic form $\om$ extends to a closed form $\Om$ on $P_\ga$, that we can arrange to be symplectic by adding the pullback of a suitable form on the base $S^2$.

The bundle $P_\ga\to S^2$ carries two canonical cohomology classes, the first Chern class $c_1^{\Ver}$ of the vertical tangent bundle and
the {\it coupling class} $u_\ga$, which is
the unique class that extends the
fiberwise
symplectic class $[\om]$ and is such that $u_\ga^{n+1}=0$.
Then, with notation as in (\ref{eq:QH}), we define
\begin{equation}\labell{eq:S}
\Ss(\ga): = \sum_{\si, i} \bla \xi_i\bra^{P_\ga}_\si\;\xi_i^{*_M}\otimes q^{-c_1^{\Ver}(\si)} 
t^{-u_{\ga}(\si)}\in H_*(M)\otimes \La,
\end{equation}
(cf. \cite[Def.~11.4.1]{MS2}.)  
Thus $\Ss(\ga)$ is obtained by \lq\lq counting" all section classes in $P_\ga$
through one fiberwise constraint $\xi_i$. 
As in Seidel's original paper~\cite{Sei},
one can also think of it as the Floer continuation map
around the loop $\ga$; cf. \cite[\S12.5]{MS2}. 
We also note for later use that 
for all $b\in H_*(M)$
\begin{equation}\labell{eq:Sa1}
\Ss(\ga)*b = \sum_{\si,i} \bla b,\xi_i\bra^{P_\ga}_\si\;\xi_i^{*_M}\otimes q^{-c_1^{\Ver}(\si)} 
t^{-u_{\ga}(\si)},
\end{equation}
i.e. it is given by counting section classes with two fiberwise constraints. 

We will also often use the fact that if $\si_0$ is 
the homology class of a section of $P_\ga\to S^2$, then every other section  in $H_2(P_\ga)$ may be written as $\si_0+\be$ for a unique $\be\in H_2(M)\subset H_2(P_\ga)$.  Note  that
\begin{equation}\labell{eq:c1v}
c_1^{\Ver}(\si_0 + \be) = c_1^{\Ver}(\si_0) + c_1(\be)
\end{equation}
where as usual  $c_1$ denotes $c_1(TM)$.

 The first parts of the following lemma are due to Seidel~\cite{Sei}.
 
 \begin{lemma}\labell{le:cN}   Let $N$ be the minimal Chern number of $M$ and $\ga\in \pi_1(\Ham\,M)$. \SSS
 
 \NI{\rm  (i)}  If $N\ge 2n+1$ then $\Ss(\ga) = \1\otimes \la$; 
 \SSS
 
 \NI {\rm  (ii)} If  $n+1\le N\le2n$ then  $\Ss(\ga^{N}) = \1\otimes \la$.
 Moreover  $\Ss(\ga) = \1\otimes \la$ provided that
  $(M,\om)$ is not strongly uniruled;\SSS
 
  \NI{\rm  (iii)}  If $N=n$ and  if $(M,\om)$ is not strongly uniruled then 
  $\Ss(\ga) = \1\otimes \la + pt\otimes q^n\mu$ for some $\la\ne 0$;\SSS
  
  \NI{\rm  (iv)} If $(M,\om)$ is negatively monotone then 
  $\Ss(\ga) = \1\otimes r_0t^{\eps_0} + x$ where $x\in \Qq_-,\; \nu(x)< \eps_0$.\SSS
  
   \NI{\rm  (v)} If $\om = 0$ on $\pi_2(M)$ then $\Ss(\ga) = (\1+x)\otimes r_0t^{\eps_0}$ where  $x\in H_{<2n}(M)[q]$.
 \end{lemma}
 \begin{proof} Suppose that
  $N \ge 2n+1$ and that
  $\bla a,b,c\bra^M_\be\ne 0$, where $a,b,c\in H_*(M)$ have degrees $2d_a,2d_b,2d_c$ respectively.
  Then $n + c_1(\be) + d_a+d_b+d_c=3n$. This equation has no solution if $|c_1(\be)|\ge 2n+1$. Hence $c_1(\be)=0$, and
  $d_a+d_b+d_c=2n$.  Moreover, if $a=pt$ 
  then $d_b=d_c = n$ which is impossible because 
  $\bla pt,M,M\bra^M_\be=0$ when $\be\ne 0$.
  Hence $M$ is not strongly uniruled, and so by Lemma~\ref{le:Qq}
   all units in $QH_*(M)$ have the form
  $\1\otimes \la+x$ where $\la\ne0,x\in \Qq_-$.
  Therefore to prove (i) it suffices to show that $x=0$.
  
 Next observe that the Seidel element is given by  invariants of the form $ \bla \xi \bra^{P_\ga}_\si$ where $\xi\in H_*(M)$. 
If $\si$ is a section class, then
 $c_1(TP)(\si) = c_1^{\Ver}(\si) + 2$, 
 where the $2$ is the Chern class of the tangent bundle to the section. Hence
 the dimension over $\C$  of the moduli space of parametrized $\si$-curves in $P_\ga$   is $n+1 + c_1^{\Ver}(\si) + 2$. Therefore
  the above invariant can be nonzero only if 
\begin{equation}\labell{eq:degxi}
2c_1^{\Ver}(\si) + \deg \xi=0.
\end{equation}
In particular, we must have $-n\le c_1^{\Ver}(\si) \le 0$.

Now suppose that the coefficient of $\1$ in $\Ss(\ga)$ is nonzero.  Then there is a nonzero invariant of the form $ \bla pt \bra^{P_\ga}_{\si_0}$, which implies that the corresponding section $\si_0$ has $c_1^{\Ver}(\si_0) = 0$.
Equation (\ref{eq:c1v}) then implies that
every other section $\si_0+\be$ either has $c_1^{\Ver}(\si_0+\be)=0$ so that this section also contributes to the coefficient of $\1$, or
is such that $c_1^{\Ver}(\si_0+\be) =  c_1(\be)$ has absolute value at least $N$.  Hence if $N>n$  equation (\ref{eq:degxi}) shows
that the class $\si_0+\be$ cannot contribute to   $\Ss(\ga)$. Thus $\Ss(\ga) = \1\otimes \la$. This proves (i).
 
 In case (ii) the argument in the previous paragraph shows 
  that all sections of $P_\ga$ contributing to 
 $\Ss(\ga)$ have the same value of  $c_1^{\Ver}$, say $-d_\ga$.
  For short let us say that $\ga$ is at level $d_\ga$.  If $d_\ga>0$ then 
  $\Ss(\ga)\in \Qq_-$.  Since $\Ss(\ga)$ is invertible,
   it follows from
   Lemma~\ref{le:Qq} that this is possible only if $(M,\om)$ is strongly uniruled.  This proves the second statement in (ii).  To complete the 
  proof of 
   (ii) we must show that  $\ga^N$ always lies at level $0$ even if $(M,\om)$ is strongly uniruled.
 
To see this, observe that for any loops $\ga,\ga'$ one can form the fibration
 corresponding to their product $\ga\,\ga'$ by taking the fiber connect sum 
 $P_\ga\# P_{\ga'}$.  Thus
 if $\si_\ga$ is a section of $P_\ga$ and $\si_{\ga'}$ is a section of 
 $P_{\ga'}$ one can form a section of $P_{\ga\ga'}$ by taking their connected sum $\si_\ga\#\si_{\ga'}$.   Under this operation the vertical Chern class adds.  Thus
 the sections that contribute to $\Ss(\ga\,\ga')$ either have level
 $d_\ga+d_{\ga'}$ or have level $d_\ga+d_{\ga'}-N$, depending on which of these numbers lies in the allowed range $[0,n]$. 
 Hence if $\ga$ lies at level $d$, for all $k\ge 1$ there is $m\ge 0$ such that
 $kd-mN\in [0,n]$.  Since $N>n$ the only possible solution of this equation when $k=N $ is $m=d$.   Therefore $\Ss(\ga^N)$ lies at level $0$.
 This proves (ii).
 
Almost the same argument proves (iii). Note that if there is a section at level $d> 0$ then all sections that contribute to $\Ss(\ga)$ lie at this level.
On the other hand if there is a section at level $0$ there might also be a section at level $n$.  It follows as before that 
$\Ss(\ga)$ must have a section at level $0$.  
Therefore $\Ss(\ga) = \1\otimes \la + pt\otimes q^n\mu.$  Moreover,
our assumption on $QH_*(M)$ implies that  $\la\ne 0$.

To prove (iv) note that $(M,\om)$ is not uniruled, so that 
by  Lemma~\ref{le:Qq}
$\Ss(\ga) = \1\otimes \la + x$ where $\la \ne 0$. 
Moreover, if $\si_0$
is a section class of $P$ of minimal energy with 
$\bla pt\bra_\si\ne 0$ then every other section class $\si$
with $\bla a\bra_\si\ne 0$ for some $a\in H_*(M)$
 has  the form $\si=\si_0+\be$ where $c_1(\be)\le0$. 
 If $c_1(\be)=0$ then $\om(\be)=0$ and $a=pt$. These invariants contribute to the coefficient of $t^{\eps_0}$ in $\la$, where $\eps_0 = -u_{\ga}(\si_0)$.   On the other hand if $c_1(\be)<0$ then $\om(\be)>0$ and hence these invariants contribute to terms in $x$ with valuation $<\eps_0-\om(\be)$.

Finally, if $\om$ vanishes on $\pi_2(M)$, then, because 
all Gromov--Witten invariants vanish,  the quantum multiplication is undeformed. Moreover, all sections of $P_\ga$ have the same energy. Hence 
$\Ss(\ga)= (\1 + x)\otimes r_0t^{\eps_0}$ where $0\ne r_0\in \Q$ and
$x\in \Qq_-$.  This proves (v).
 \end{proof}

\begin{defn}\labell{def:si}   If $\Ss(\ga)=\1\otimes \la+x$  we define $\si_0$
to be the section class of $P_\ga$ of minimal energy that contributes 
nontrivially to the coefficient $\la$.  Moreover we write
$\la = \sum_{i\ge 0}r_i\,t^{\eps_i}$
 where $r_0\ne 0$ and
$\eps_i>\eps_{i+1}$ for all $i$.
\end{defn}

Thus 
\begin{equation}\labell{eq:r0}
r_0 = \bla pt \bra^{P_\ga}_{\si_0},\qquad
\nu(\la) = \eps_0 =-u_\ga(\si_0).
\end{equation}
Note that  $c_1^{\Ver}(\si_0)=0$ by equation (\ref{eq:degxi}).
The conditions in Proposition~\ref{prop:sch} are chosen 
so that the moduli space of sections in class $\si_0$ is compact.
In what follows we shall often simplify notation by writing $P$ 
instead of $P_\ga$.
\MS

\NI{\bf Proof of Proposition~\ref{prop:sch}.}
  Suppose that $\om|_{\pi_2}=0$. Then
$\Ss(\ga) = (\1 + x)\otimes r_0t^{\eps_0}$ by
 Lemma~\ref{le:cN} (v). 
We will show that
$\eps_0 = 0$.  The claimed result then follows from 
Proposition~\ref{prop:des} (i:d). 

Since $\om|_{\pi_2(M)} = 0$ the spaces of sections 
that give the coefficients of $\Ss(\ga)$ are compact
(when 
 parametrized as sections) since there is no bubbling.
 To be more precise, note that it suffices to consider
 almost complex structures $J$ on $P: = P_\ga$ that are compatible with the 
 fibration $\pi:P\to S^2$ in the sense that $\pi$ is holomorphic and the restriction $J_z$ of $J$ to the fiber over $z\in S^2$ is
 $\om$-compatible for each $z\in S^2$.  Then all stable maps in the
 compactification of a space of sections consist of a section together with some bubbles in the fibers.  But if $\om|_{\pi_2(M)} = 0$ there are
 no $J_z$-holomorphic spheres for any $z$. 
 Hence in this case  the sections in class $\si_0$ form a {\it compact} 
 $2n$-dimensional manifold $\Mm: = \Mm(\si_0)$. (Here we assume that the elements of $\Mm$ are parametrized as sections.) Therefore
  there is a commutative diagram
 $$
 \begin{array}{ccc} \;\; S^2\times \Mm
 &\stackrel{ev}\longrightarrow& \;\;P\\
pr_1 \downarrow &&\pi \downarrow\\
\;\;S^2 & = & \;\;S^2\end{array}
$$
where the top arrow is the evaluation map.  Moreover, $\si_0
= ev_*[S^2\times \{pt\}]$.  

Now observe that the Gromov--Witten invariant
$r_0=\bla pt \bra^{P}_{\si_0}$ of equation (\ref{eq:r0})
is just the degree\footnote
{In fact, since $r_0$ must be an integer it must be $\pm 1$.  For its product with the corresponding integer for $\ga^{-1}$ must be $1$.}
 of the map
$\{z\}\times \Mm\to M=\pi^{-1}(z)$ that is induced by $ev$. 
Thus our hypotheses imply that $ev$ has nonzero degree.  Since the coupling class $u_\ga$ satisfies $u_\ga|_M=[\om]$ and 
$u_\ga^{n+1}=0$, it 
 easily follows that there is a class $a\in H^2(\Mm)$ such that 
$ev^*(u_\ga) = pr_{2}^*(a)$. 
Hence 
$$
\eps_0=\nu(\la) = - 
u_\ga(\si_0)  = -\int_{S^2\times \{pt\}} ev^*u_\ga =  0
$$
as required.   This proves (i).

Next consider  case (ii). 
Lemma~\ref{le:cN} implies that in all cases considered here 
 $\Ss(\ga) = \1\otimes \la + x$
where $\nu(x)<\nu(\la)$.  Hence the spectral invariants descend provided that  $\nu(\la) = 0$.  But for generic $\pi$-compatible
 $J$ on $P$
we claim that the moduli space $\Mm$ is again compact.   Hence the previous argument  applies to show that $\nu(\la) = 0$.

 To prove the claim,
note that the only possible stable maps in the compactification $\oMm$ of $\Mm$ consist of the union of a section in the
class $\si_0-k\be$ where $k>0$ together with some bubbles in classes $k_i\be$.  Suppose first that $c_1(M)=\ka [\om]$ where $\ka<0$.  Then 
because $\om(\be)>0$, $c_1(\be)\le -N$.   But if $J$ is generic it induces a generic $2$-parameter family $J_z,z\in S^2,$ of almost 
complex structures on $M$.  Therefore such a class $\be$ could be represented only if $-N\ge c_1(\be)  \ge 2-n$, which is possible only if $N\le n-2$.  Since we assume $N\ge n-1$
there are no bubbles in this case.  On the other hand if 
$c_1(M)=\ka [\om]$ on $\pi_2(M)$ where
$\ka> 0$ then 
 $c_1(\be) \ge N$.  Hence  $c_1(\si_0-k\be)\le 2-N$.  
But the section is embedded and hence is regular for generic $J$.
Thus  for such a section to exist
we must have  $n+1 + (2-N) \ge 3$, i.e. $n\ge N$. (This is precisely the argument in \cite{EP1}.) 
This proves (ii).

If $N\ge n+1$  but $(M,\om)$ is strongly uniruled, then we can apply the above argument together with Corollary~\ref{cor:des}
to $\ga^N$ to conclude that the 
asymptotic spectral 
invariants descend.  
\QED

\subsection{Calculating the coupling class.}\labell{s:out}

An essential ingredient of Schwarz's argument
 is that the vanishing of $\om$ on $\pi_2$ implies that there is no bubbling
 so that  the moduli space 
 $\Mm$ is a compact manifold. 
One cannot replace $ S^2\times \Mm$ here
 by the universal curve $\oMm_{0,1}(\si_0)$ over some compactified virtual moduli cycle $\oMm(\si_0)$
since the \lq\lq bundle" $f:\oMm_{0,1}(\si_0)\to \oMm(\si_0)$ is singular over the higher strata of $\oMm(\si_0)$.  The  argument
we give below  shows
 that if the relevant Gromov--Witten invariants of $M$ and $P: = P_\ga$ vanish then this potential twisting does not effect $ev^*(u_\ga)$ too much, so that this class still has zero integral over the fiber of $f$.  Observe that it is not enough here that the invariants of $M$ vanish; we need some control 
 on the  invariants of $P$ in either section or fiber classes.

Our reasoning
is very similar to that in \cite[\S3]{Mcq}; see also \cite{Mcv}. There we were investigating conditions under which the ring $H^*(P)$ 
splits as a product, i.e. it is isomorphic as a ring to
 $H^*(M)\otimes H^*(S^2)$.  As the next lemma makes clear, what we need here is a partial splitting of this ring.

\begin{lemma}\labell{le:H2}  Suppose 
$\Ss(\ga) = \1\otimes \la+x$ and that
there is an 
element $H\in H_{2n}(P_\ga;\R)$ such that\SSS

\NI {\rm (a)}  $H\cap [M]$ is Poincar\' e dual in $M$ to $[\om]$;\SSS

\NI {\rm (b)} $H\cdot \si_0=0$;\SSS

\NI {\rm (c)} $H^{n+1} = 0$.\SSS

\NI Then $\nu(\la) = 0$.
\end{lemma}

\begin{proof}  Let $u: = {\PD}_P (H)$, the Poincar\' e dual in $P: = P_\ga$ to the divisor class $H$. Then (a) implies that $u|_M = [\om]$, while (c) implies that $u^{n+1} = 0$.  Hence $u$ is the coupling class
$u_\ga$.   Therefore
$$
\nu(\la): = -u_\ga(\si_0) = -u(\si_0) = -H\cdot \si_0  = 0,
$$
as required.
\end{proof}

It is possible that an element $H\in H_{2n}(P)$ with the above properties exists whenever $N\ge n+1$ and $\Ss(\ga)=\1\otimes \la$.  However,  we can only prove this under an additional assumption\footnote
{
This condition, though essential to the proof, is very technical 
and there seems no intrinsic reason why it should be necessary.}
 on the structure of $QH_*(M)$ that we now explain.  Roughly speaking we require that the divisor classes in $H_*(M)$ carry all the nontrivial quantum products.  To be more precise, we make the following definition.

\begin{defn}\labell{def:D}
Let $\Dd$ be the subring of $H_{ev}(M)$ generated by the divisors, i.e. the
  elements in $H_{2n-2}(M)$.  We say that $QH_*(M)$ {\bf 
  satisfies condition (D)}
  if  $\Dd$ has an additive complement $V$ in $H_{ev}(M)$ such that the following conditions hold
  for all $d\in \Dd, v\in V$ and $\be \in H_2(M)$
  $$
  (a)\qquad d\cdot v=0;\qquad \mbox{and }\;\;\;(b)\qquad\bla d,v\bra^M_\be = 0.
  $$
  \end{defn}
 This condition is used in Lemma~\ref{le:DH} in order to allow certain computations to be done recursively.
 
 \begin{rmk}\rm (i) 
  This condition is trivially satisfied if either the quantum product is undeformed (since then  $\bla d,v\bra^M_\be=0$ always)
  or if $\Dd=H_{ev}(M)$.   
  \SSS
  
  \NI (ii) When $N\ge n+1$ the only classes $\be$ 
  for which $\bla d,v\bra^M_\be $ could be nonzero have
$c_1(\be)=0$.  (If $N=n+1$ there is one potentially nonzero invariant with two insertions and $c_1(\be)\ne 0$, namely $\bla pt,pt\bra^M_\be$.  But $pt\in \Dd$ and  so this does not affect 
condition (b) in Definition~\ref{def:D}.)
\end{rmk}

    By condition  (a) in Definition~\ref{def:D}
     we may choose a basis $\xi_i,0\le i\le m,$ for $H_{ev}(M)$ so that
  the first $m_1$ elements span the subring $\Dd$ while the others span $V$.  
  Hence the first $m_1$ elements of the dual basis $\{\xi_{i}^{*_M}\}$ will also span $\Dd$.
In other words, for all  $d\in \Dd$ and  $\be \in H_2(M)$
 \begin{equation}\labell{eq:DV}
    \bla d,\xi_i\bra^M_\be\ne 0\; \Longrightarrow\; \xi_i^{*_M}\in \Dd\quad\mbox{ for all } i.
\end{equation}
Notice also that because there is an open subset of $H_{2n-2}(M)$ consisting of elements $D$ such that $D^n\ne 0$ we may assume that $\xi_0=\1$ and that
  for $1\le i\le m_1$ each $\xi_i = D^k$ for some such $D$.  Then each
  of the corresponding dual elements 
  $\xi_i^*$ is also a sum of elements of the form $D^k$.

Now consider the map $
s: H_*(M)\to H_{*+2}(P)
$   defined by the identity
\begin{equation}\labell{eq:s}
s(a)\cdot_P v: = {\textstyle \frac 1{{r_0}}}\, \bla a,  v\bra^P_{\si_0},
\qquad v\in H_*(P),
\end{equation}
where $r_0$ is as in equation (\ref{eq:r0}).

\begin{prop}\labell{prop:H}  Suppose that $\Ss(\ga) = \1\otimes \la + x$, let  $\si_0$ be as in Definition \ref{def:si}  and define $s$ as above.
Then the element $H: = s(\PD_M(\om))$ 
satisfies the 
 conditions of Lemma~\ref{le:H2} in each of the following situations:\SSS
 
 \NI
 {\rm (i)}   $N\ge n+1$,  $\Ss(\ga) = \1\otimes \la$ and 
$QH_*(M)$ 
  satisifies condition  {\rm (D)};\SSS
 
 \NI {\rm (ii)}  $N=n$,  $(M,\om)$ is not strongly uniruled, $\rank\, H_2(M)>1$, and condition (D) holds;\SSS
 
  \NI {\rm (iii)}  $(M,\om)$ is negatively  monotone and condition {\rm (D)} holds;\SSS
  
  \NI {\rm (iv)}   $QH_*(M)$ is undeformed and $(M,\om)$ is not spherically monotone with $\rank\, H_2(M) = 1$.  
\end{prop}

We defer the proof to \S\ref{s:GW}.
 
\begin{cor}
Theorems~\ref{thm:main1} and ~\ref{thm:main2} hold.
\end{cor}
\begin{proof} 
Theorem~\ref{thm:main1} concerns the case when $QH_*(M)$ is undeformed. Then $\Ss(\ga)$ always has the form $\1\otimes \la' + x$ where $x$ is nilpotent, since these are the only invertible elements in $QH_*(M)$.  Moreover, if $N\ge n+1$, $\Ss(\ga) = \1\otimes \la$ by Lemma \ref{le:cN} (ii).  Hence we may suppose that
 the conditions of either part (i) or part (iv) of Proposition~\ref{prop:H} hold.  Therefore this result, together with Lemma~\ref{le:H2} and 
Corollary \ref{cor:des}, prove the theorem.   Note that if $n\le N$   we cannot conclude 
 that the spectral invariants descend, since it is possible that
 $\nu(x)>0$ for some $\ga$.

Similarly, part (iv) of 
Theorem~\ref{thm:main2} follows from part (iii) of
 Proposition~\ref{prop:H}.
Next suppose  that $N\ge n+1$ and, if $N\le 2n$,
 that $(M,\om)$ is not strongly uniruled.  Then 
$\Ss(\ga) = \1\otimes \la$ by Lemma~\ref{le:cN}.  Moreover
$\nu(\la)=0$ by
Lemma~\ref{le:H2}.  Hence the spectral invariants descend
by part (i) of Proposition \ref{prop:des}. This proves part (i) of 
Theorem~\ref{thm:main2}.  Parts (ii) and (iii) of this theorem follow similarly. In case (ii) we have $\Ss(\ga^N) = \1\otimes \la$  while in case (iii) 
$\Ss(\ga) = \1\otimes \la+pt\, \otimes q^n\mu$.  In either case 
 $\nu(\la)=0$.
Therefore the asymptotic  spectral invariants descend by 
Corollary~\ref{cor:des}.  
\end{proof}

\section{Calculations of Gromov--Witten invariants}\labell{s:GW}

This section contains the proof of Proposition~\ref{prop:H}.

\subsection{Preliminaries.}\labell{s:gen}

We shall use the following identity of 
Lee--Pandharipande~\cite{LP}:
 \begin{eqnarray}
 \labell{eq:LPii}
ev_i^*(H)& = &ev_j^*(H) + (\al\cdot H)\,\psi_j -\sum_{\al_1+\al_2=\al}
(\al_1\cdot H)\, D_{i,\al_1\,|\,j,\al_2},
\end{eqnarray}
where $\psi_i$ is the first Chern class of the cotangent bundle to the domain at the $i$th marked point and 
$H\in H_{2n}(P)\cong H^2(P)$ is any divisor.\footnote
{
Here  $\al\cdot H$ denotes the intersection number of  two homology 
classes $\al, H$ that lie in complementary degrees.  If $u\in H_i(P)$
where $i$ is arbitrary, we shall denote by $Hu$ the cap product $H\cap u\in H_{i-2}(P)$.
Further, when $H=M$ is the class of the fiber, we shall write $u\cap M$ for the cap product when considered as an element in $H_{i-2}(M)$.
This last distinction is not very important since $H_*(M)$ injects into $H_*(P)$ by the result of \cite{Mcq}.}
  Lee--Pandharipande
were working in the algebraic context and hence 
interpreted both sides as elements of an appropriate Picard group.
Thus
provided that we are working with stable maps that have at least three marked points,
 $D_{i,\al_1\,|\,j,\al_2}$ is the divisor consisting of all stable maps with two parts, one in class $\al_1$ and containing $z_i$ and the other 
 in class $\al_2 = \al-\al_1$ containing  $z_j$.
 
We shall interpret  (\ref{eq:LPii})  as an 
 identity for  Gromov--Witten invariants.   Thus taking $i,j=1,2$
 this equation
states that for all classes $u,v,w\in H_*(P)$ and  all 
$\al\in H_2(P)$,
\begin{eqnarray}\labell{eq:LP3}
\bla H u, v, w\bra_{\al} &=& \bla u, Hv, w\bra_{\al}
+(\al\cdot H)\;\bla u, \tau v, w\bra_{\al}\\ \notag
&& \qquad - 
\sum_{j, \al_1+\al_2=\al} ({\al_1}\cdot H)\;
\bla u,\eta_j,\dots\bra_{\al_1}\;\bla \eta_j^*, v, \dots\bra_{\al_2},
\end{eqnarray}
where the sum is over the elements of a basis $\{\eta_j\}$ of $H_*(P)$ (with dual basis $\{\eta_j^*\}$) and all decompositions $\al_1+\al_2=\al$ of $\al$, and where the dots indicate that the constraint $w$ may be in either factor
except if $\al_2=0$ in which case it must be in the second
 factor for reasons of stability. Note that $\tau$ here denotes a descendent invariant.  In fact we shall only use (\ref{eq:LPii}) in cases 
 when $H\cdot\al=0$ so that this term vanishes.  If it does not, one should 
 get rid of the $\tau$ insertion by using the identity $\psi_i = D_{i|jk}$,
 where $D_{i|jk}$ denotes the divisor consisting of all stable maps with two parts, one containing the point $z_i$ and the other containing the points $z_j,z_k$.
  (But often one then gets no information since the term on LHS appears in the expansion for $\tau$.) For further discussion of this point as well as 
 a brief proof of (\ref{eq:LPii}) in the symplectic context
 see \cite{Mcu}.

Notice that if we apply an identity such as (\ref{eq:LP3}) 
to even dimensional classes $a_i$ we only need to consider 
even dimensional $\eta_j, \eta_j^*$. We will make this restriction now in order to avoid irrelevant  considerations of sign.
Also, we will simplify the arguments that follow by choosing a basis 
$\{\eta_j\}
$ for $H_{ev}(P)$ of special form.   
We start with a basis $\xi_i, i\in I,$ for $H_{ev}(M;\R)$ that satisfies the condition in (\ref{eq:DV})
and extend this to a basis for $H_{ev}(P;\R)$ by adding elements 
$\xi_j^*, j\in I,$
so that for all $i,j$
\begin{equation}\labell{eq:basis}
\xi_i\cdot_P \xi_j^* = \de_{ij},\quad 
\xi_i^*\cdot_P \xi_j^*=0.
\end{equation}
Thus $\xi_i$ is a fiber class but $\xi_j^*$ is not. 
 Note that this basis
$\{\xi_i, \xi^*_j\}$ for $H_{ev}(P)$ is self-dual.  Further,
$\xi_i^{*M}=\xi_i^*\cap M$.

With this choice of basis the sum in (\ref{eq:LP3}) breaks into two,
depending on which of $\eta_j,\eta_j^*$ is a fiber class.  To analyze the resulting product terms, we
frequently use the following fact about invariants of 
$P$ in a fiber class $\be$.

\begin{lemma}\labell{le:fib}  
Suppose that $a,b\in H_*(M)$, $v,w\in H_*(P)$ and $\be\in H_2(M)$.  Then:\SSS

\NI {\rm (i)}
$
\bla a,b,v\bra^P_\be = 0.$\SSS
 
 \NI{\rm (ii)} 
 $
 \bla a,v,w\bra^P_\be =\bla a,v\cap M,w\cap M\bra^M_\be.
 $
 \end{lemma}
\begin{proof}   This is part of \cite[Prop~1.6]{Mcq}.
These statements hold because, as is shown in \cite{Mcq},
 one can calculate these invariants using an almost complex structure and perturbing $1$-form that
are compatible with the fibration  $\pi:P\to S^2$.
 Hence every element in the virtual moduli cycle is represented by a curve that lies in a single fiber. (i) is then immediate, since $a,b$ can be represented in different fibers.  Similarly, (ii) holds because 
 every $\be$-curve through $a$ must lie in the fiber containing $a$.
 The fact that $v\cap M, w\cap M$ do not vanish means that
 these cycles take care of the needed transversality normal to this fiber.
\end{proof}

To make our argument work we also need information about certain section invariants of $P$.  When $N\ge n+1$ the following lemma suffices; the proof is easy since it is based on a dimension count.
 
\begin{lemma}\labell{le:easy}  Suppose  that $N\ge n$ and that
if $N=n$ then $(M,\om)$ is not strongly uniruled.  Suppose further that
$\Ss(\ga) = \1\otimes \la + pt\otimes q^n \mu$, and let  
 $\si=\si_0-\be$ where $\om(\be)>0$.  Then for all
 $a\in H_{<2n}(M)$:
\MS

\NI {\rm (i)} For all $w\in H_{*}(P)$, 
$\bla a,w,M\bra_{\si}=0$   unless $c_1^{\Ver}(\si)=0$ and
$\deg w+\deg a=2n$;\SSS

\NI {\rm (ii)}
 $\bla  a,b,M\bra_{\si}=0$   for all  $b\in H_*(M)$;
\SSS

\NI {\rm (iii)} For all $w\in H_{*}(P)$, 
$\bla a,w,M\bra_{\si}$
depends only on $w\cap M$.
\end{lemma}
\begin{proof} Statement (i) holds by a dimension count as in the proof of Lemma~\ref{le:cN}.  (ii) holds because an invariant of this form
with only two nontrivial fiber insertions is determined by $\Ss$; namely
 by (\ref{eq:Sa1}) $\bla a,\xi_i\bra_{\si}$ is the coefficient of 
 $\xi_i^{*_M}\otimes q^{-c_1^{\Ver}(\si)}\,t^{-u_\ga(\si)}$ in
$$
\Ss(\ga)*a = (\1\otimes \la + pt\otimes q^n\mu)*a = a\otimes \la.
$$
Since ${-u_\ga(\si)}> \nu(\la)$ this vanishes. Note that for this
argument to apply
we need  that either $\mu=0$ (which happens when $N\ge n+1$) or  $pt*a=0$, i.e. $M$ not strongly uniruled. 

Statement (iii) is an immediate consequence of (ii) 
because  any two classes  $w,w'$  with
 $w\cap M = w'\cap M$  differ by $w-w'\in H_*(M)$.\end{proof}

\subsection{Technical lemmas.}\labell{s:tech}

This section contains two rather technical results about the 
vanishing of certain Gromov--Witten invariants needed for some of the cases considered by Proposition~\ref{prop:H}.  However, they are not needed when $N\ge n+1$, and the reader might
do well to read the next section first, coming back to 
this section later.

If $N=n$ we also need the following lemma about the fiber invariants of $P$.

\begin{lemma}\labell{le:hard}  Suppose  that $N=n$, that
$(M,\om)$ is not strongly uniruled, and that $\rank\, H_2(M)>1$. Suppose further that
$\Ss(\ga) = \1\otimes \la + pt\otimes q^n\mu$ where $\la,\mu\in 
\La^{univ}$.    Then for any  classes $s_1,s_2\in H_2(P)$ and any $\be\in H_2(M)$,
$\bla  s_1,s_2\bra_\be=0.$
\end{lemma}
\begin{proof}  A dimension count shows that the invariant is zero unless $c_1(\be) = n$.  If there are any nonzero invariants of this form,
 choose $\be$ with minimal energy
$\om(\be)$ so that $\bla  s_1,s_2\bra_\be\ne 0$.

If $s_1,s_2$ are both fiber classes then the invariant vanishes by
Lemma~\ref{le:fib}.  If just one (say $s_1$) is a fiber class then the same lemma implies that the invariant equals $\bla  s_1,pt\bra^M_\be$ which vanishes because $M$ is not strongly uniruled.  Hence
 $\bla  s_1,s_2\bra_\be$  does not depend on the choice of section classes $s_i$.   
 
 We now make a specific choice of $s_1$.
 Since $\rank\, H_2(M)>1$ there is $a\in H^2(M)$ such that $a(\be) = 0$.  Let $F\in H_{2n}(P)$ be any extension of $\PD_M(a)\in H_{2n-2}(M)$.
(Since $H_*(M)$ injects into $H_*(P)$ by \cite{Mcq} such a class exists.)
Choose $b\in H_2(M)$ such that $F\cdot_P b=pt$, and let 
$v\in H_4(P)$ be any extension of $b$.  Then $vF: = v\cap F\in H_2(P)$  is a section class of $P$.  Therefore if $\si$ is any section class, it suffices to show that
$$
\bla  vF, \si, H\bra_\be=0,
$$
where $H\in H_{2n}(P)$ is chosen so that $H\cdot \be = 1$.

To prove this, apply (\ref{eq:LPii}) with $i=1, j=2$.  We obtain
\begin{eqnarray*}
\bla  vF, \si, H\bra_\be &=& \bla  v, F\si, H\bra_\be + (\be\cdot F)\,
\bla  v, \tau\si, H\bra_\be \\
&&\qquad - \;\sum_{i,\be_1+\be_2=\be} (\be_1\cdot F)\, \Bigl(\bla  v, \xi_i,\dots\bra_{\be_1}
\bla  \xi_i^*, \si, \dots\bra_{\be_2}\\
&&\qquad\qquad\qquad + \; \bla  v, \xi_i^*,\dots\bra_{\be_1}
\bla  \xi_i, \si, \dots\bra_{\be_2}\Bigr),
\end{eqnarray*}
where the dots indicate that the $H$ insertion could be in either factor.
The first term in RHS vanishes because $F\si$ is a multiple of a point.
 Hence Lemma~\ref{le:fib} implies that
 $ \bla  v, F\si, H\bra_\be= \bla  v\cap M, pt, H\cap M\bra^M_\be$
 which vanishes
because  $(M,\om)$ is not uniruled.  The second term on RHS
vanishes since $\be\cdot F=0$.
 Further in the sum neither $\be_1$ nor $\be_2$ is zero because of the factor $\be_1\cdot F$.  Since $c_1(\be)=n=N$ one of the $\be_i$ has $c_1=0$ and the other has $c_1=n$. (Other possibilities such as $c_1(\be_i) = -n$ can be ruled out by dimensional considerations
  as in the proof of Proposition~\ref{prop:sch}.)

Consider the first sum.
Since $\xi_i$ is a fiber class  and $\dim (v\cap M) = 2$ there are two possibilities; either $c_1(\be_1) = 0$ and $\dim \xi_i = 2n$ or
$c_1(\be_1) = n$ and $\dim \xi_i = 0$.  In the first case, $\xi_i^*$ is a section class so that $\bla  \xi_i^*, \si, H\bra_{\be_2}=0$ by the minimality of $\om(\be)$.  Therefore such a term does not contribute. 
But the second case also does not contribute because $(M,\om)$ is not strongly uniruled.  Therefore the first sum vanishes.

Now consider the second sum. Applying Lemma~\ref{le:fib} again, we find that because $\xi_i$ is a fiber class the second factor here equals
$\bla  \xi_i, pt , H\cap M\bra^M_{\be_2}$.  Hence it vanishes 
by hypothesis on $M$.
\end{proof}

When all $3$-point Gromov--Witten invariants in $M$ vanish the same argument gives 
the following conclusion.
 
\begin{lemma}\labell{le:hard1}  Suppose that the quantum 
multiplication in $M$ is undeformed and that $\rank\, H_2(M) > 1$.
Then for all nonzero $\be\in H_2(M)$, $v\in   H_*(P)$ and
 $F\in H_{2n}(P)$ 
we have:\SSS

\NI {\rm (i)}  $\bla F^k,v\bra_\be=0$; 

\NI{\rm (ii)} $\bla \si,v \bra_\be=0$ for all section classes $\si$.
\end{lemma}
\begin{proof}  Consider (i).
If any of $F^k,v$ are fiber classes then the invariant reduces to an invariant in $M$ and hence vanishes by assumption.
Therefore, as above we may assume that $F\cdot\be = 0$.
Now suppose that there is some nonzero invariant as in (i)
 and choose a minimal $k$ such that 
$\bla F^k,v\bra_\be\ne0$.  Then $k>1$ by the divisor axiom. 
Now choose $H$ so that $\be\cdot H = 1$ and
expand  $\bla F^k,v,H\bra_\be$ as before.
The first term vanishes by the minimality of $k$ and the second
since $\be\cdot F=0$.  Moreover, in the sum neither of $\be_1,\be_2$  
vanish.  Therefore, because each term in the sum has a fiber constraint in at least one of its factors the sum vanishes.  This proves (i).

To prove (ii) let $F$ be any extension of the Poincar\'e dual to $[\om]$.  Then $F^n\cap M=pt$. Hence $F^n=\si$ is a section class.
Therefore $\bla \si,v \bra_\be=\bla F^n,v \bra_\be=0$ by (i).
\end{proof}

\subsection{Proof of Proposition~\ref{prop:H}.}\labell{s:pr}

  For simplicity we will now assume  that $[\om]$ is normalized so 
that $\int_M\om^n = 1$. Further, we set
 $h: = \PD_M(\om)$ so that $h^n=pt$.\footnote
 {
 Since $[\om]\in H^2(M)$ need not be a rational class, $h\in H_{2n-2}(M;\R)$. Therefore, in this section one should assume that homology groups have coefficients $\R$ unless otherwise indicated.
}
  Recall that we define $s:H_*(M)\to H_{*+2}(P)$ by:
$$
s(a)\cdot_P v: = {\textstyle \frac 1{{r_0}}}\, \bla a,  v,M\bra^P_{\si_0},\qquad v\in H_*(P).
$$
 In particular, $s(M)\cdot_P pt = \frac 1{r_0}\, \bla M,  pt,M\bra^P_{\si_0}=1$ so that $s(M)=P$.

\begin {lemma}\labell{le:s} Suppose that  
$\Ss(\ga)=
\1\otimes \la + x$ where for all $b\in H_{<2n}(M;\R)$ either 
$\nu(x*b)< \nu(\la)$  or $x*b = q y$ where 
$y\in QH_*(M)$ involves only nonnegative powers  of $q$.  Then:\SSS

\NI
  {\rm (i)} $s(pt) = \si_0$,  \SSS

\NI {\rm (ii)}  
  $s(a)\cap M = a$ for all $a\in H_*(M)$.
 \end{lemma}
\begin{proof}
 $\si_0$ was chosen so that   $\bla pt,M\bra^P_{\si_0} = 
 \bla pt\bra^P_{\si_0}=r_0\ne 0$.
 A count of dimensions shows that
$\bla pt,  v\bra^P_{\si_0}\ne 0$ only if $v$ is a divisor class.
Thus for all $v\in H_{2n}(P)$ the divisor axiom implies that
$$
s(pt)\cdot  v \;= \; 
{\textstyle \frac 1{{r_0}}}\, \bla pt,  v\bra^P_{\si_0} \;= \;
{\textstyle \frac 1{{r_0}}} ({\si_0}\cdot v) \,\bla pt\bra^P_{\si_0} 
\;= \;{\si_0}\cdot v.
$$
This proves (i).

Since we have already checked the case $a=[M]$ in (ii), it suffices to 
take $a\in H_{<2n}(M)$.  
Observe   that by equation (\ref{eq:Sa1})  
 $\bla a,\xi_i\bra_{\si_0}$
 is the coefficient of 
$\xi_{i}^{*_M}\otimes 
t^{\eps_0}$ in
$$
\Ss(\ga)*a \;= \; (\1\otimes \la +x)*a \;= \; a\otimes \la + x*a,
$$
The hypothesis on $x$ implies that $x*a$ does not contribute to this coefficient. 
Hence
$
\bla a,\xi_i\bra_{\si_0}=
r_0\,a\cdot_M\xi_i,
$
and so 
$$
{r_0} \,(s(a)\cap M)\cdot_M \xi_i\; =\; {r_0}\, s(a)\cdot_P \xi_i \;\stackrel{def}=\;
\bla a,\xi_i\bra_{\si_0} \;= \;{r_0}\, a\cdot_M\xi_i
$$
where the second equality holds by the definition (\ref{eq:s}) of $s$.
Hence
$
s(a)\cap M = a$ as required.
\end{proof}

There are many different cases in Proposition~\ref{prop:H}.  In an attempt to avoid confusion we will first prove parts (i) and (ii). Thus we assume that $N\ge n$, with some further conditions when $N=n$.

  \begin{lemma}\labell{le:less}  Suppose that the hypotheses of Lemma~\ref{le:easy} hold, and that if  $N=n$ those of
  Lemma~\ref{le:hard} hold as well. 
Then   $\bla h,\si_0,M\bra_{\si_0}=0$.
\end{lemma}
\begin{proof} 
Given any 
divisor class in $P$  extending $h$ we may add a suitable multiple of $M$ to obtain a class $K\in H_{2n}(P)$ such that
$$
K\cap M = h,\quad K\cdot\si_0 = 0.
$$
We then find by formula (\ref{eq:LPii}) that
\begin{eqnarray*}
\bla h,\si_0,M\bra_{\si_0} &= &
\bla KM,\si_0,M\bra_{\si_0}\\
&=& \bla M,K\si_0,M\bra_{\si_0}+(\si_0\cdot K)\bla M,\tau\si_0,M\bra_{\si_0}\\
&&\quad -\; \sum \bigl((\si_0-\al)\cdot K\bigr)\; \Bigl(
\bla M,\xi_i,\dots\bra_{\si_0-\al}
\bla \xi_i^*,\si_0,\dots\bra_{\al}\\
&&\qquad\qquad\qquad\quad
+ \;\bla M,\xi_i^*,\dots\bra_{\si_0-\al}
\bla \xi_i,\si_0,\dots\bra_{\al}\Bigr),
\end{eqnarray*}
where the dots indicate that the $M$ insertion could be in either factor.
Note that the first two terms vanish because $K\cdot\si_0=0$.
  Suppose there is a nonzero contribution from the first sum.  Because the first factor has two fiber constraints, Lemma~\ref{le:fib} implies that
  $\si_0-\al$ must be a section class. Moreover  
   $$
0\ne    (\si_0-\al)\cdot K = -\al\cdot_M h = \om(\al).
$$
Thus  $\om(\al)> 0$ since there is a nonzero $\al$-invariant.  But then 
$$
 \bla M,\xi_i\bra_{\si_0-\al}=  \bla M,\xi_i,M\bra_{\si_0-\al}
$$ vanishes by 
Lemma~\ref{le:easy}(ii) except possibly if $\dim \xi_i = 2n$.  
But this case can occur only if $N=n$ and $c_1(\al) =n$.  But then $\xi_i^*$ is a section class, so that the second factor $\bla \xi_i^*,\si_0,\dots\bra_{\al}$ vanishes by Lemma~\ref{le:hard}.  
   Therefore in all cases this sum makes no contribution.

Now consider the second sum.  Suppose first that $\be: = \si_0-\al$ were a fiber class. Note that $\be\ne 0$ because  the sum is multiplied by
 $\be\cdot K$.
Then, because $M$ is a fiber constraint, Lemma~\ref{le:fib} implies that the 
$\be$ invariant  is either $\bla M,\xi_i^*\cap M, M\cap M\bra_{\be}$,
or is
$\bla M,\xi_i^*\cap M\bra_{\be}$, which both vanish when $\be\ne 0$
because the first marked point is not constrained.  Thus the only nonzero terms have $\si_0-\al$ a section class, and hence have the form
$\bla M,\xi_i^*,M\bra_{\si_0-\al}^P
\bla \xi_i,pt\bra_{\al}^M$ where  $\om(\al)
\ne 0$ since $K\cdot\si_0=0$.    Again there are two cases.
If $N>n$ then 
Lemma~\ref{le:easy}(i) implies that 
this can be nonzero only if $\deg \xi_i^* = 0$.
But $\xi_i^*\ne pt$ since it is not a fiber class. Therefore this is impossible.
On the other hand if $N=n$ then the first factor might be nonzero, but the second has to vanish since $M$ is not uniruled. 
Therefore in all cases the second sum vanishes as well.
 \end{proof}

  \begin{cor}\labell{cor:Hsi} Under the conditions of the above lemma,
  if $H: = s(h)$, we have $H\cdot \si_0 = 0.$
  \end{cor}
  \begin{proof}  By the definition of $s$ in (\ref{eq:s}),
  $s(h)\cdot \si_0 = \bla h,\si_0,M\bra_{\si_0}=0$.
  \end{proof}

  \begin{lemma}\labell{le:DH}
    Suppose that $(M,\om), \Ss(\ga)$ and $N$ satisfy the hypotheses of   Lemma~\ref{le:less} and that  $QH_*(M)$ satisfies condition {\rm (D)}.  
 Then:
\SSS

\NI  
  {\rm (i)} 
  For all section classes $\si=\si_0-\be$ where $\om(\be)>0$
  and all $k$ we have
  $\bla F^k,a,M\bra_\si=0$ whenever 
   $F\in H_{2n}(P)$ and $a\in H_{<2n}(M)$. \SSS

  \NI {\rm (ii)} $\bla H^{n+1-k},h^k, M\bra_{\si_0}=0$ for all $k$.
  \end{lemma}
  \begin{proof} Consider (i).  
   Choose $\si$ of minimal energy (i.e. $\om(\be)$ is maximal)
   and then the minimal $k$ so that   
   $\bla F^k,a,M\bra_\si\ne0$ for some $a$ and $F$.  Note that  $k>1$
   since 
   $\bla F,a,M\bra_\si= (F\cdot \si) \bla a,M\bra_\si$  vanishes by Lemma~\ref{le:easy}(ii).  

Again by Lemma~\ref{le:easy}(ii) we may add an arbitrary fiber class to $F^k$ without changing $\bla F^k,a,M\bra_\si$.  Hence, by replacing $F$ by $F-cM$ for suitable $c$, we may arrange
 that $F\cdot\si=0$.   
 Now use (\ref{eq:LPii}) as in the previous lemma, moving $F$ from the first constraint to the second.  Because $F\cdot\si=0$ we find as before that
\begin{eqnarray*}
 \bla F^k,a,M\bra_{\si} &=&  \bla F^{k-1},Fa,M\bra_{\si}-
 \sum (\al\cdot F) \Bigl(\bla F^{k-1},\xi_i,\dots\bra_{\al}
\bla \xi_i^*,a,\dots\bra_{\si-\al}\\
&&\qquad\qquad
+ \;\bla F^{k-1},\xi_i^*,\dots\bra_{\al}
\bla \xi_i,a,\dots\bra_{\si-\al}\Bigr).
\end{eqnarray*}
The first term is zero by our choice of $k$.  Consider the first sum. If  $\al$ is a section class, then because of the factor $\al\cdot F$ we may assume that
$\al\ne \si$.   But then $\al$ has less energy than $\si$ since $\om(\si-\al)>0$ so that  
the $\al$ invariant vanishes  by the minimality of the energy of $\si$.  So
we may suppose that  $\al$ is a fiber class, in which case the other factor is an invariant 
  $\bla \xi_i^*,a,M\bra_{\si-\al}$  in a section 
  class with smaller energy than $\si$.
  We claim that condition (D) implies that the product of these factors must  still vanish.
  
For suppose not.  Then  by Lemma~\ref{le:fib}
  $$
  \bla F^{k-1},\xi_i\bra^P_\al = 
    \bla f^{k-1},\xi_i\bra^M_\al\ne 0
    $$
 where $f: = F\cap M$.  Therefore   (\ref{eq:DV}) implies that
     $\xi_i^{*_M} = \xi_i^*\cap M\in \Dd$.
Thus $\xi_i^*\cap M$ is a linear combination of elements of the form 
$g^k, g\in H_{2n-2}(M)$.  But if we write $g=G\cap M$  for some 
$G\in H_{2n}(P)$
all invariants of the form 
$\bla G^k,a,M\bra_{\si-\al}$ vanish by the minimality condition on the 
energy of $\si$.   But $\xi_i^*\cap M = G^k\cap M$.  Hence
   $\bla \xi_i^*,a,M\bra_{\si-\al}=0$   by Lemma~\ref{le:easy}(iii).
  
   Thus the first sum vanishes.
   Now consider the second sum.  The second invariant has at least two fiber constraints. Hence $\si-\al$ must be a section class, so that the invariant vanishes
  by Lemma~\ref{le:easy}(ii).  This proves (i).
  
  Now consider (ii).  When $k=n$ the invariant is
  $\bla H, pt,M\bra_{\si_0}$ which vanishes because $H\cdot\si_0=0$.
  Suppose that it does not vanish for all $k$ and
 choose  the maximal $k$ 
  for which it is nonzero.  Expand
  $\bla H^{n+1-k},h^k, M\bra_{\si_0}$ by transferring one $H$ to the second constraint. As usual the first two terms in the expansion vanish and we obtain
\begin{eqnarray*}
\bla H^{n+1-k},h^k, M\bra_{\si_0} &=&  -
 \sum (\al\cdot H) \Bigl(\bla H^{n-k},\xi_i,\dots\bra_{\al}
\bla \xi_i^*,h^k,\dots\bra_{\si_0-\al}\\
&&\qquad\qquad
+ \;\bla H^{n-k},\xi_i^*,\dots\bra_{\al}
\bla \xi_i,h^k,\dots\bra_{\si_0-\al}\Bigr).
\end{eqnarray*}
Consider the first sum.  If $\al$ is a section class then we may suppose $\al\ne \si_0$ because of the factor $\al\cdot H$.  
Therefore it has smaller energy than $\si_0$ so that the $\al$-invariant vanishes by part (i) of this lemma.  On the other hand if $\al$ is a fiber class then as above condition (D) implies that $\xi_i\in \Dd$ 
and the second factor is a sum of terms of the form
$\bla K^{n-k+1}, h^k,M\bra^P_{\si}$ where $\si$ has less energy than $\si_0$.  Therefore this factor vanishes by part (i).  Therefore the first sum vanishes.  
But the second factor in the second sum has at least  two fiber 
constraints. Therefore $\si_0-\al$ is a section class, and its  energy is
less than that of $\si_0$ since $\om(\al)\ne 0$.  Hence   this factor must vanish by Lemma~\ref{le:easy}(ii).
Thus the RHS of the above expression vanishes.  Therefore the
LHS is zero also, contrary to hypothesis.   The result follows.  \end{proof}
  
  \begin{cor}\labell{cor:DH}
If  the hypotheses of Lemma~\ref{le:DH} hold then $H^{n+1} = 0$.
\end{cor}
\begin{proof} 
Putting $k=1$ into claim (ii)  of Lemma~\ref{le:DH}
we find that
$\bla h,H^{n},M\bra^P_{\si_0}=0$.  But by equation  (\ref{eq:s}) 
this  is a multiple of the intersection $s(h)\cdot H^n$.  Since $s(h) = H$, we obtain $H^{n+1}=0$.
\end{proof}

\NI {\bf Proof of Proposition~\ref{prop:H}.}\,\,
If $N\ge n+1$, condition (D) holds and $\Ss(\ga) = \1\otimes \la$,
then the hypotheses of Lemma~\ref{le:easy} hold.  Hence we may apply Lemmas~\ref{le:less} and \ref{le:DH}.  Therefore, the conditions of Lemma~\ref{le:H2} hold by
Corollaries \ref{cor:Hsi} and \ref{cor:DH}  and  Lemma~\ref{le:s}.  This proves (i).
To prove (ii) note that the extra conditions here precisely match the conditions of Lemma~\ref{le:hard}.  Hence the proof goes through
as before.

Now consider part (iii) of the proposition.
The assumption here is that $(M,\om)$ is negatively monotone and that condition (D) holds. We saw in Lemma~\ref{le:cN} that in this case
 $\Ss(\ga) = \1 \otimes r t^{\eps_0} + x$ where $\nu(x)<\eps_0.$ Hence the conclusions of Lemma~\ref{le:s} hold.
  
Further Lemma~\ref{le:easy}(ii) holds (though part (i) may not). 
To see this, recall that $\bla a,\xi_i^{*M},M\bra_{\si}$ 
is the coefficient of 
$\xi_i \otimes q^{-c_1{\Ver}(\si)}\,t^{-u_\ga(\si)}$ 
in $\Ss(\ga)*a$. Because $\nu(x)<\nu(\la) $ in $\Ss(\ga)$, 
this must vanish when $-u_\ga(\si)>-u_\ga(\si_0)$.
Therefore if we write $\si = \si_0-\be$, the invariant vanishes when
 $\om(\be)>0$.  Therefore part (iii) of this lemma also holds.

It is now easy to check that the proof of
 Lemma~\ref{le:less} goes through.  The argument for the vanishing of the first sum needs no change (note that $c_1(\al)<0$ here so there is no exceptional case); 
 while that for the second sum works using the fact that $(M,\om)$ is not uniruled instead of Lemma~\ref{le:easy}(i).
It remains to check the proof of Lemma~\ref{le:DH}. But this holds 
as before, provided that condition (D) hold.
This completes the proof of part (iii).

Finally consider part (iv).
 The assumption here is that $QH_*(M)$ is undeformed, and if 
 $\rank\,H_2(M) =1$ that $(M,\om)$ is not positively monotone. 
 Since condition (D) holds when $QH_*(M)$ is undeformed, the
 latter case follows from (i) if $c_1=0$ and from (iii) otherwise.  Hence we may assume that
 $\rank\,H_2(M) >1$.  Therefore  the conditions of 
 Lemma~\ref{le:hard1} hold. Further,  
because  $\Ss(\ga) = \1\otimes \la + x$ has degree $2n$,
 all terms in $x$ have a positive coefficient of $q$. When
 $x*b=x\cap b$ this remains true for $x*b$.  
 Hence the conclusions of Lemma~\ref{le:s} hold.

Next we claim that 
Lemmas~\ref{le:less} and \ref{le:DH} hold.
To see this, we go though the proofs of these lemmas
 using Lemma~\ref{le:hard1} instead of Lemma~\ref{le:easy} to show that the requisite terms vanish. Note for example that
the fiber invariants in the  expansion in Lemma~\ref{le:less}
contain factors of the form $\bla \xi^*,\si\bra_\al$ which vanish by Lemma~\ref{le:hard1}(ii). Similarly, in
Lemma~\ref{le:DH} we may use Lemma~\ref{le:hard1}(i).  We do not need condition (D) because  in a product such as
$$
\bla H^{n-k},\xi_i,M\bra_{\al}
\bla \xi_i^*,h^k\bra_{\si_0-\al}
$$
with $\al\in H_2(M)$  the first factor vanishes, and so 
we do not need to worry about the second factor.
 Hence the proof goes through  as before.\QED

\section{Examples}\labell{ss:ex}

We now prove Proposition~\ref{prop:bl}. We shall calculate in 
the subring
$QH_{2n}(M) \cong H_*(M)\otimes \La^{univ}$, i.e. by fixing the degree
of the elements considered we can forget the coefficients $q^i$.

Recall the following construction from K\c edra~\cite{K2} 
and \cite{Mcbl}.
Let $(X,\om_X)$ be a symplectic manifold.
For each map $\al:S^2\to X$, let $gr_\al$ be the graph $(z,\al(z))\in S^2\times X$ of $\al$, and let $\tau$ be an area form on $S^2$ of area $1$.
Choose a constant $\mu_0$ so that 
\begin{equation}\labell{eq:mu0}
\Om: = \mu_0\, pr_1^*(\tau) + pr_2^*(\om)
\end{equation}
is nondegenerate on $gr_\al$.
Denote by $(\TP,\TOm_\de)$ the $\de$-blow up of 
$(S^2\times X,\Om)$ along $gr_\al$.  (The parameter 
$\de$ refers to the symplectic area of a line in the exceptional divisor.)
Then the symplectic bundle $
\pi:\TP\to S^2$ has fiber $M: = (\TX,\Tom_\de)$ and corresponds to  
a Hamiltonian loop $\ga_\al\in \pi_1(\Ham(M))$.

\begin{lemma}\labell{le:bl}  Let $(X,\om_X)$ be a symplectic manifold of dimension $2n\ge 4$.
Given a map $\al:S^2\to X$ define
$\mu_\al: = \int_{S^2}\al^*\om_X$ and  $\ell_\al: = c_1^X(\al)$.
Suppose that $\ell_\al\ge 0$ and that
  at least one of  $\mu_\al,\ell_\al$ is nonzero.
Then
 the Seidel element $\Ss(\ga_\al)\in QH_{2n}(M)$ of the loop $\ga_\al$ defined above has 
the form $\1\otimes \la + x$ where $\nu(\la)\ne 0$
for small nonzero $\de$.
\end{lemma}
\begin{proof}  Let $D$ be the exceptional divisor in $\TP$.  
Denote the trivial section of $S^2\times X$ by $\si_1: = S^2\times \{p\}$.
If $p\notin gr_\al$ this lifts to a section $\Tsi_1$ of $\TP\to S^2$. If $\eps$ denotes the class of the line in the fibers of the exceptional divisor, then every section class may be written as
$\Tsi_1 - m\eps + \be$ where $\be \in H_2^S(X)$.  If $m=0$ then this class is pulled back from $S^2\times X$, and   
because
 $\TP$ is obtained from $S^2\times X$ by blowing up along a
(complex) curve with nonnegative Chern class, we may apply 
the results of Hu \cite{Hu}.  Thus  for all $a\in H_*(X)$
\begin{equation}\labell{eq:Hu}
\bla a\bra^{\TP}_{\Tsi_1+\be} = \bla a\bra^{S^2\times X}_{\si_1+\be}.
\end{equation}
Hence $\bla a\bra^{\TP}_{\Tsi_1+\be}=1$ if $\be = 0$ and vanishes otherwise.  

Define $\ka: = -u_\ga(\Tsi_1)$.  The above argument shows that
the coefficient $\la$ of $\1$ in
$\Ss(\ga)$ contains the term $t^{\ka}$.
There might be other classes  
$\Tsi_1 - m\eps + \be$ that contribute to the coefficient of $\1$ in
$\Ss(\ga)$ but these all  appear
with the coefficient
$t^{\ka +m\de -\om(\be)}$ where $m\ne 0$.
Thus $\nu(\la) $ either equals $\ka$ or equals  $\ka +m\de -\om(\be)$
for some $m\ne 0$ and $\be\in H_2(X)$.

We now  calculate $\ka$.  The coupling class $u_\ga$ has the form $\TOm_\de +c_\de\, pr_1^*\tau$ where $c_\de$ is chosen so that
$\int_{\TP} u_{\ga}^{n+1}=0$.    Further, because 
$\TOm_\de|_{\Tsi_X} = \Om|_{\si_X}$ by construction and $\int_{\si_1}\Om=\mu_0$ by equation (\ref{eq:mu0}),
$$
-\ka := u_{\ga}(\Tsi_1) = c_\de + \mu_0.
 $$
We showed in \cite{Mcbl} (see also \cite{Mcsn}) that if $V: = \frac 1{n!}\int_X\om^n$
then
$$
{\textstyle \frac 1{(n+1)!}\int_{\TP} (\TOm_\de)^{n_1} = \mu_0(V - v_\de)-
v_\de\left(\mu_\al -\frac {\ell_\al}{n+1} \de\right)},
$$
where $v_\de: = \frac{\de^n}{n!}$ is the volume of the ball
 cut out of $X$.
But
$$
{\textstyle 0 = \frac 1{(n+1)!}\int_{\TP}(\TOm_\de +c_\de\, pr_1^*\tau)^{n+1} =
\left(\frac 1{(n+1)!}\int_{\TP} (\TOm_\de)^{n+1}\right) + c_\de \,(V-v_\de),}
$$
since $V-v_\de$ is the volume of $(M,\Tom_\de)$.
Therefore
$$
{\textstyle  -\ka: = c_\de +\mu_0 =  \frac{v_\de}{V-v_\de}
\left(\mu_\al -\frac {\ell_\al}{n+1} \de\right).}
$$
Thus,  provided that at least one of
$\mu_\al,\ell_\al$ are nonzero,  $\ka$ is a
rational function of $\de$ with isolated zeros.   Therefore
neither $\ka$ nor $\ka + m\de -\om(\be)$ vanish for 
sufficiently small $\de\ne0$.  
 The result follows.
\end{proof}

\NI
{\bf Proof of Proposition~\ref{prop:bl}.}\,\, 
The previous lemma proves (i),
and so it remains to consider the case when $X$ has dimension $4$
and $[\om]$ and $c_1$ vanish on $\pi_2(X)$. 
It is shown in 
 \cite[Prop~6.4]{Mcbl}
 that under the given hypotheses on $X$
  every Hamiltonian bundle $(M,\Tom_\de)\to (\TP,\TOm)\to S^2$ is constructed by blowing up some section $\si_X$ of some Hamiltonian bundle $(X,\om)\to (P,\Om)\to S^2$. 
 
 Let us denote by 
 $\Tga_X$ the loop in $\pi_1(\Ham\,M)$ corresponding to the $M$-bundle $\TP\to S^2$ and by $\ga_X$ that corresponding to its blow down $P\to S^2$.
As in the proof of Proposition~\ref{prop:sch},  
$\Ss(\ga_X) = r_0 \1$ for all 
 $\ga_X\in \pi_1(\Ham\, X)$.  Therefore there is at least one\footnote
{There may be several such since each coefficient of $\Ss(\Tga)$ 
is a sum of contributions from all sections with given values of $c_1^{\Ver}$ and $\mu_{\Tga}$.}  
section
$\si_{X0}: = \si_X+\be$ of $P$ such that 
$\bla pt\bra^{P}_{\si_{X0}} \ne 0$.   
Moreover both $c_1^{\Ver}$ and the coupling class $u_X$ of $X$ vanish on the section $\si_{X0}$, and hence on all other sections of $P$, in particular on $\si_X$.  

The section classes in $\TP$ have the form $\Tsi_{X0} - m\eps+\be$, where $\Tsi_{X0}$ is the lift of $\si_{X0}$ and $\be\in H_2^S(M)$.  For such a class to contribute to 
  the corresponding Seidel element $\Ss(\Tga)$ we need 
  $-2\le c_1^{\Ver}(\Tsi_X - m\eps + \be) = -m\le 0$.
Moreover, since we may choose $\Om$ so that
the section $\si_X$ of $(P,\Om)$ is symplectic, 
Hu's results\footnote
{
We do not need to use Hu \cite{Hu} here. All that matters is that
the coefficient of $\1$ in $\Ss(\Tga_X)$ is nonzero, which follows from the fact that $M$ is not uniruled.} imply that
$$
\bla pt\bra^{\TP}_{\Tsi_{X0}}= \bla pt\bra^{P}_{\si_{X0}} \ne 0.
$$ 
Therefore if $\ka: = -u_{\Tga_{X}}(\Tsi_{X0})$
   $$
 \Ss(\Tga_X) = r_0\1\otimes t^{\ka} + \xi\otimes  t^{\ka+\de} + r\, pt\otimes  t^{\ka+2\de},\quad \mbox{ where }\xi\in H_2(M).
 $$
We can now repeat the calculation of $\ka$ in  Lemma~\ref{le:bl}.
All we need change is  the interpretation of the constant $\mu_0$,
which we now {\em define} to be the area $\Vol(P,\Om)/\Vol(X,\om_X)$
of the fibration $P\to S^2$.
Hence $\ka = 0$.

Now observe that as in \cite[\S2]{Mcu}
 $$
 \eps*\eps = -pt + \eps \otimes  t^{-\de},
 $$
  and that quantum multiplication  (in $M$) by 
  the elements of $H_{\le 2}(X)\subset H_{\le 2}(M)$ is undeformed. 
  In particular, $pt*\eps = 0$ so that
 $$
 \nu\big((\eps\otimes t^\de)^k\bigr) = \de,\qquad k\ge 1.
 $$ 
 Hence,  if we decompose the class $\xi$ appearing in $\Ss(\Tga_X)$
 as $\xi = s\eps + \xi'$ for some $s\in \Q, \xi'\in H_2(X)$,
we have  $\xi'*\eps = \xi'*pt = 0$, and
 we easily find that
 $$
 0\le \nu(\Ss(\Tga^k))\le 2\de,\qquad k\ge 1.
 $$
Hence the asymptotic invariants descend by Proposition~\ref{prop:des}. 
\QED

\begin{rmk}\rm (i) It is perhaps worth pointing out that the equality in (\ref{eq:Hu}) does not always hold if you blow up along 
the graph of a class with
$c_1$ negative.  For example, suppose that you take $\al:=L$ 
to be the line in $X: = \C P^2$.  Then it is not hard to check that
$\Ss(\ga_{\al}) = \1\otimes t^\ka +(L-E) \otimes t^{\ka+\de}$, where $E$ denotes the exceptional divisor in $M$. (The second term comes from counting sections in class $\Tsi_0 - E$.)  Let us normalize the symplectic form on $X$ so that $\om_X(L)=1$. Since $\ga_{\al}^{-1} = \ga_{-\al}$, we must then have
$$
\Ss(\ga_{-\al}) = \Ss(\ga_{\al})^{-1} =(- E + pt\otimes t^{\de})
t^{1-2\de-\ka}.
$$
Thus $\ga_{\al}$, which is formed by blowing up along the graph of $-L$,
has a Seidel element  in which  
the coefficient of $\1$ vanishes.
These calculations are carried out in detail in \cite[\S5]{Mcv}.
Observe that $\ga_{-\al}$ is three times the 
generator of $\pi_1(\Ham M)$ that is called $\al$ in  
\cite{Mcv}. Thus the element of $QH_*(M)$ called 
$Q^{-3}$ in  \cite{Mcv} has the form $\Ss(\ga_{\al}) t^{\ka'}$ for appropriate $\ka'$. See also \cite[Remark~1.8]{Mcbl}.\SSS

\NI (ii) One should be able to use the methods of  \cite{Mcu}
and Lai~\cite{Lai} to show that in the situation of Lemma~\ref{le:bl}
 classes with $m\ne 0$ do not contribute to $\Ss(\ga)$.  This calculation will be carried out elsewhere.
\end{rmk}

\end{document}